\documentclass[12pt]{article}

\usepackage{amssymb,amsmath}
\usepackage{amsthm}
\usepackage{epsfig,graphicx}

\newcommand{\dis}{\displaystyle}

\newcommand{\be}{\begin{equation}}
\newcommand{\ee}{\end{equation}}
\newcommand{\bea}{\begin{eqnarray}}
\newcommand{\eea}{\end{eqnarray}}
\newcommand{\bean}{\begin{eqnarray*}}
\newcommand{\eean}{\end{eqnarray*}}
\newcommand{\ba}{\begin{array}}
\newcommand{\ea}{\end{array}}

\newcommand{\Rr}{{(\mathbb R)}}
\newcommand{\C}{{\mathbb C}}

\newcommand{\Z}{{\mathbb Z}}

\renewcommand{\phi}{\varphi}

\newcommand{\DC}{\Delta^{\Omega^c}_D}
\newcommand{\RC}{R^{\Omega^c}_D}
\newcommand{\R}{R^{\Omega}_D}
\newcommand{\DB}{\Delta^{\Omega^c}_B}
\newcommand{\RB}{R^{\Omega^c}_B}
\newcommand{\DL}{{\mathcal D}\ell(k)}
\newcommand{\SL}{{\mathcal S}\ell(k)}

\newtheorem{theo}{{\sc Theorem}}[section]
\newtheorem{cor}[theo]{{\sc Corollary}}
\newtheorem{lem}[theo]{{\sc Lemma}}
\newtheorem{prop}[theo]{{\sc Proposition}}

\newtheorem{rem}[theo]{{\sc Remark}}

\title{Inverse resonance problem  for
\boldmath{$\Z_2$}-symmetric\\ analytic obstacles
in the plane\thanks{%
Research partially supported by  NSF grant \#DMS-0071358.}}

\author{Steve Zelditch\thanks{%
Department of Mathematics,
Johns Hopkins University,
Baltimore, MD 21218, USA
(zelditch@math.jhu.edu).}}

\begin{document}

\maketitle

\begin{abstract}
We given an exposition of a proof that a mirror symmetric
configuration of two convex analytic obstacles in $\Rr^2$ is
determined by its Dirichlet resonance poles. It is the analogue
for exterior domains of the proof that a mirror symmetric bounded
simply connected analytic plane domain is determined by its
Dirichlet eigenvalues.  The proof uses 'interior/exterior duality'
to simplify the argument.
\end{abstract}

\section{Introduction} This article is part of a
developing series \cite{Z1, Z2} concerned with the inverse
spectral problems for analytic plane domains.  It is essentially
the  lecture we presented at the IMA workshop on Inverse Spectral
Problems in July, 2001 and represents the state of our knowledge,
methods and results at that time. We apply these methods to prove
an analogue for exterior domains of the result proved in \cite{Z2}
for interior ones, namely  that a mirror symmetric configuration
of two convex analytic obstacles in $\Rr^2$ is determined by its
Dirichlet or Neumann resonance poles.  The proof combines the
known result that wave invariants of an exterior domain are
resonance invariants (see \S \ref{POISSON}) with  the method of
\cite{Z1, Z2} for calculating the wave invariants explicitly in
terms of the boundary defining function. In keeping with the
expository nature of the lecture, we give a detailed exposition of
the background results on the Poisson relation for exterior
domains and of  the main steps in \cite{Z1, Z2} in the calculation
of wave invariants (\S \ref{COEFFS}) and on determining   the
domain from its wave invariants (\S \ref{DOMAIN}).

 The motivating problem is   whether
analytic domains are determined by their spectra for Dirichlet or
Neumann boundary conditions. This inverse spectral problem can be
posed for both  interior and exterior domains. For the interior
problem, we assume the domain is a bounded  simply connected plane
domain  and ask whether its Dirichlet or Neumann spectrum
determines the boundary up to rigid motion. For the exterior
problem, we assume the domain is the complement of two bounded,
simply connected obstacles and ask whether we can determine the
pair of obstacles from the resonance poles (or scattering phase)
of the exterior Dirichlet or Neumann Laplacian. Since the method
is similar for both boundary conditions, we assume the boundary
conditions are Dirichlet in the exterior domain.

In the previous paper \cite{Z2}, the interior inverse spectral problem is
studied for analytic domains with a certain mirror symmetry
$\sigma.$ The symmetry is assumed  to fix a bouncing ball orbit (as a set)
 and to reverse its orientation; i.e. if the domain is translated and rotated so that
the bouncing ball orbit lies along the y-axis with its midpoint at
the origin, then $\sigma(x, y) = (x, - y).$  In addition,  we
assume
 the length of the bouncing ball orbit is a fixed number $L$ and that the orbit
satisfies a non-degeneracy condition. The main result (see Theorem
(\ref{ONESYMINT}))  is that the interior domain satisfying some
generic conditions is determined by its Dirichlet spectrum among
other such domains. We note that this result is stronger than the
one in \cite{Z3} in having eliminated one of the two symmetries.
We should emphasize that the proof does not immediately imply the
result for the other mirror symmetry $\sigma(x, y) = (-x, y)$,
which fixes $\gamma$ pointwise and preserves the orientation.

Our main purpose in this paper is to extend   the result to
  connected exterior domains of the form
\begin{equation} \label{OBST} \Omega^c := \Rr^2 - \Omega,\;\; \mbox{with}\;\;
\Omega = \{{\mathcal O} \cup \tau_{x, L} {\mathcal O}\} \end{equation}
where ${\mathcal O}$ is a  bounded,  simply connected
analytic domain  and where
$$\tau_{x, L} = \; \mbox{reflection through the line}\; \langle \nu_x, x + 1/2 L \nu_x - y \rangle = 0.$$
Here,  $x \in {\mathcal O}, $  $\nu_x$ is the outward unit normal
to ${\mathcal O}$ at $x$ and $L > 0$ is a given positive number.
Thus, the obstacle consists of two (non-intersecting) isometric
components which are mirror images of each other across a common
orthogonal segment of length $L$.  Such an obstacle is called a
(symmetric) 2-component scatterer. The segment is the projection
to $\Rr^2$ of a bouncing ball orbit $\gamma$ of the exterior
billiard problem. Throughout, the notation $\Omega$ refers to a
bounded domain and $\Omega^c$ is the unbounded complement.

\begin{figure}[htb]%f1\vglue-15pt
\centerline{\includegraphics{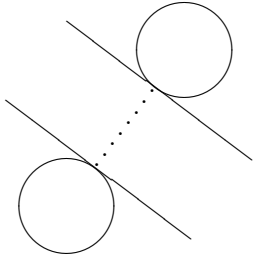}}\vglue-15pt
\caption{Mirror symmetric 2-obstacle scatterer.}
\end{figure}

The problem is thus to recover the obstacle ${\mathcal O} \cup \tau_{x, L} {\mathcal O}$
 from the  set of resonance poles $\{ \lambda_j\}$  of the exterior Dirichlet Laplacian  $\DC$, i.e. the
poles of the analytic continuation of its resolvent $$\RC(k + i \tau) = (\DC + (k + i \tau)^2)^{-1}$$
to the logarithmic plane.
 To be precise,
we consider the class ${\mathcal OBSTA-}$\break
${\mathcal CLE}_{\Z_2, L}$ of
two-component obstacles   satisfying:

\eject

\begin{itemize}

\item  (i)   ${\mathcal O}$
is simply\; connected and real  analytic;

\item (ii) $ \gamma$ is a non-degenerate bouncing
ball orbit, whose length $L$ is isolated in $ Lsp(\Omega) \cup
Lsp(\Omega^c)$.

\end{itemize}

 We denote by ${\mathcal RES}
(\Omega^c)$ the set of resonance poles  of the Laplacian $\DC$  of
the domain $\Omega^c $ with Dirichlet boundary conditions.

We now state the main results.

\begin{theo} \label{ONESYM}    ${\mathcal RES}: OBSTACLE_{\Z_2, L} \mapsto \C^{{\bf N}}$ is 1-1.

\end{theo}

In fact, in combination with a result of M. Zworski \cite{Zw}, the
proof shows more: the obstacle is determined by the resonances
close to the real axis which are associated to the bouncing ball
orbit $\gamma$  of the exterior billiard problem. We note that
assumption (ii) is stronger than the one in the interior case (cf.
Theorem (\ref{ONESYMINT}) in demanding multiplicity one for $L$ in
the combined interior and exterior length spectra. By using the
original proof in \cite{Z1, Z2}, one could remove $Lsp(\Omega)$
from  the assumption and just demand multiplicity one in
$Lsp(\Omega^c)$. We do not do so, because we wish to present a
simpler proof than the one in \cite{Z1, Z2}, as we will explain
below in the introduction.

\begin{figure}[htb]%f2\vglue9pt
\centerline{\includegraphics{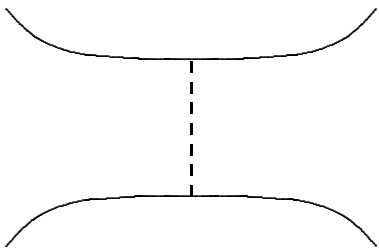}} \caption{Bouncing ball
orbit.}
\end{figure}

As a corollary, we prove that the obstacle is determined by
its scattering phase.
We denote by $ S_{D}(\lambda)$ the Dirichlet scattering operator for $\Omega^c$ and by
$ s_D(\lambda) = \det S_{D}(\lambda)$  the scattering phase.

\begin{cor}  If two exterior domains $\Omega$ in ${\mathcal OBSTACLE}_{\Z_2, L}$ have the same scattering phase
$s_D(\lambda)$, then they are isometric.   \end{cor}

We note that the length $L$ of the segment between the components
has been marked, since apriori it is not a resonance invariant.
Under some additional assumptions, it is  a resonance invariant
and one does not have to mark it. For instance, if one assumes
that ${\mathcal O}$ is convex, then the bouncing ball orbit
$\gamma$ is the unique periodic reflecting ray of the exterior
domain up to iterates, so its length is a resonance invariant.
Thus, let ${\mathcal COBSTACLE}_{\Z_2}$ denote the convex
component obstacles  in ${\mathcal OBSTACLE}_{\Z_2}$. We then have:

\begin{cor}   ${\mathcal RES}: {\mathcal COBSTACLE}_{\Z_2} \mapsto \C^{{\bf N}}$ is 1-1.  \end{cor}

For the sake of completeness, we also recall  the precise
statement of the interior result of \cite{Z2}. Let Spec$(\Omega)$
denote the spectrum of the interior Laplacian $\Delta_B^{\Omega}$
of the domain $\Omega$ with either Dirichlet or Neumann boundary
conditions.

\enlargethispage{5pt}
\begin{theo} \label{ONESYMINT}  Let
$\mathcal D_{ L}$ denote the class of  simply connected
real-analytic plane domains $\Omega$ satisfying:
\begin{itemize}

\item  (i) There exists an  isometric involution $ \sigma$ of $\Omega$;

\item (ii) $\sigma$ `reverses' a  non-degenerate bouncing ball orbit $  \gamma \to \gamma^{-1}$;

\item (iii) The lengths $2 r L$ of  all iterates $\gamma^r$ have  multiplicity one in $Lsp(\Omega)$,  and the eigenvalues of the linear
Poincare map $P_{\gamma}$ are  not roots of unity;

\end{itemize}
Then:  Spec$: {\mathcal D}_{1, L} \mapsto \Rr_+^{{\bf N}}$ is 1-1.
\end{theo}

To be cautious, we should point out that these results are part of
a work-in-progress which has not yet reached its final form. They
have been independently verified when there is an additional
mirror symmetry (i.e. each obstacle is itself `left-right' mirror
symmetric), and in this case the relevant calculations of wave
invariants are easy. When there is only one mirror symmetry, the
calculations become somewhat
 messy (see \cite{Z2})   and the proof therefore becomes rather `unstable';
at  this time of writing, some  details have been independently
verified (by R. Bacher \cite{B}),  but some have not.

Let us briefly describe the proof and the organization of the
paper. As  emphasized by Zworski, there is nothing really new in
Theorem \ref{ONESYM} beyond the result of \cite{Z1, Z2} in the
interior case (i.e. Theorem \ref{ONESYMINT}). We are just
combining the known fact that wave trace invariants of the
exterior domain are resonance invariants with the calculation of
the wave trace invariants in terms of the germ of the defining
function of $\partial \Omega$ at the endpoints of the bouncing
ball orbit. This calculation was done for the germ of  any kind of
bouncing ball orbit in \cite{Z2}, developing a method originating
in the work of Balian-Bloch \cite{BB1, BB2}.  However, it seems to
us worthwhile to collect all the facts one needs for the proof of
Theorem (\ref{ONESYM}) in one place and to explain the main steps
without all the details in \cite{Z1, Z2}.

Moreover, our exposition has one novel point:  in \S \ref{INTEXT},
we use the so-called interior/exterior duality
 to simplify formula for  the resolvent trace in \cite{Z1}.
Interior/exterior duality is  a distributional trace formula,
which we informally write for $\tau >> 0$ as
\begin{equation} \label{IO}
\begin{array}{l}
\dis Tr_{\Rr^2} [R_{ D}^{\Omega^c}(k + i \tau)
 \oplus  R_{N }^{\Omega} (k + i \tau)  - R_{0 }(k + i \tau)]\\[6pt]
\hspace{45pt}\dis =  \frac{d}{d k } \log \det \bigg(I + N(k + i
\tau)\bigg) .\end{array}
\end{equation}
Here, the determinant is the Fredholm determinant, $R_{N}^X,
R_D^X$ denote Neumann (resp. Dirichlet) resolvents on $X$ and
 the notation $Tr_X$ indicates the space on which the trace is
 taken. We write $L^2(\Rr^2) = L^2(\Omega) \oplus L^2(\Omega^c)$
 and view
 $R_{ D}^{\Omega^c}(k + i \tau)
 \oplus  R_{N }^{\Omega} (k + i \tau) $ as an operator on this
 space. For notational simplicity, we do not put in the explicit projections $1_{\Omega}, 1_{\Omega^c}$
 (i.e. the characteristic functions). Also, $N(k + i \tau)$ is a boundary integral operator which will
be defined in \S 4.
 The correct statement and proof of (\ref{IO}) will be given in  Proposition (\ref{MAIN}).
It is essentially a differentiated version of the  formula
\begin{equation} \label{EP} \frac{d}{dE} \log \det S_{ D} (E) +
\pi \frac{d}{dE} N_N (E) = \Im \frac{d}{dE} \log \zeta_{D, N} (E +
i 0) \end{equation} of Eckmann-Pillet (cf. \cite{ EP, EP2}),
relating the scattering phase of the exterior Dirichlet problem,
the eigenvalue counting function of the interior Neumann problem
and a spectral zeta function for the integral operator $N$ along
the boundary.

This reduction to the boundary trace on the right side of
(\ref{IO}) simplifies some of the technical details of \cite{Z1}
and  allows one   to deduce the inverse spectral result
simultaneously for the interior and exterior problems.
 In addition, it brings our
calculations into closer contact with the physics literature,
where the Balian-Bloch approach is now almost always applied to
the boundary trace (see e.g. \cite{AG, THS, THS2, GP}. The price
we pay is that the combination of inside and outside  requires us
to make the additional multiplicity free assumption on the
combined inner and outer length spectra to obtain the inverse
results in Theorem \ref{ONESYM}.  As mentioned above, this
assumption could be eliminated; but we feel that the
simplification in the proof is worth the extra assumption.  To
simplify the exposition, we also assume in the last step that
${\mathcal O}$ is convex. Hence, we only complete the proof for
hyperbolic bouncing ball orbits. For the general case, we refer
to~\cite{Z2}.

Let us give a brief outline of the  method  of
 \cite{Z1, Z2} and this paper  for  determining  an analytic
domain from its  spectrum.  It is based on the use of certain
spectral (i.e. resonance) invariants known as the wave trace
invariants at the bouncing ball orbit. As the name implies, the
wave trace invariants at a periodic reflecting ray $\gamma$  are
coefficients of the singularity expansion of the relative trace of
the wave group $E_D^{\Omega^c}(t) = \cos t \sqrt{\DC}$  at $t =
L_{\gamma}$ (the length). The inverse results are proved by
 explicitly
calculating   the wave trace invariants at a bouncing ball orbit
 in terms of the Taylor coefficients of the boundary defining
function at the endpoints, and then determining the boundary from
these Taylor coefficients.

Rather than studying the wave trace per se in \cite{Z1, Z2} we
follow \cite{BB1, BB2} in studying an essentially equivalent
expansion involving the resolvent. Potential theory gives an exact
formula  for the Dirichlet (or Neumann) resolvent of a domain in
terms of the free resolvent of $\Rr^2$ and of the layer potentials
associated to the domain. Formally, one may derive an infinite
series expansion known as the `multiple reflection expansion' of
\cite{BB1, BB2}, whose $M$th term corresponds intuitively to  $M$
reflections on the boundary. The trace of each term is, again
formally, an oscillatory integral corresponding to $M$-fold
periodic reflecting rays.

The advantage  of this approach  is that the formula is exact and
does not require a microlocal paramatrix construction, which is
messy and complicated for bounded domains and therefore very hard
to use in inverse spectral theory.  Moreover, the terms of the
mutiple reflection expansion are `canonical' in that the
integrands are the same for all domains; the only difference lies
in the domain one integrates over. Once it is `legalized', it
gives an explicit and even routine algorithm for calculating wave
invariants. At this time of writing, no other approach seems to
provide a workable algorithm for doing the calculations. The
disadvantages are that the legalization has two complicated parts:
one needs to estimate the remainder in the infinite series and to
regularize the oscillatory integral defined by  the $M$ term.

Our formulae for the wave invariants come from applying stationary
phase to the regularized  traces. One of the principal results is
that, for each iterate $\gamma^r$, only the term $M = rm$ of the
multiple reflection expansion contributes important inverse
spectral data, because only this term  contains the maximum number
of derivatives of the defining function of $\Omega$ at a given
order of $k^{-j}.$ The heart of the matter is then the calculation
and analysis of the wave invariants. Since one needs wave
invariants of all orders, one seemingly runs into an infinite
jungle of complicated expressions.  To tame the jungle, we
enumerate and evaluate the terms of the wave trace invariants
 using Feyman diagram techniques. It turns out that only five
 diagrams play an important role (this has been confirmed by R.
 Bacher). The remaining issue is the calculation of the
 amplitudes and their dependence on $r$. It is this latter
 dependence that is crucial in allowing one to remove a symmetry
 of the domain.

The calculation of wave trace invariants is the same for periodic reflecting
rays of interior and exterior domains since it depends only on the germ of
the boundary at the reflection points. The material in \S 5 - \S 7
summarizes the calculation in \cite{Z1, Z2}. We hope that this guide clarifies
 the calculation and makes it easier to check the details.

We close the introduction with some remarks on open problems in
the inverse spectral problem for analytic plane domains.  The most
obvious one is whether one can eliminate the remaining symmetry
assumption. It may be that one can recover the domain from wave
invariants at one closed billiard orbit, or that one has to
combine information from several orbits.  A key problem is,   how
much of the Taylor expansion of the boundary defining function one
can recover from wave invariants at one orbit?

The method of this paper and \cite{Z2} gives detailed but somewhat
undigested information about wave invariants in terms of  the
Taylor coefficients of the boundary defining function at endpoints
of a bouncing ball orbit. Another approach, suggested originally
by Colin de Verdiere \cite{CdV}, is to focus on the Birkhoff
normal form of the billiard map at the closed orbit. As was proved
by V. Guillemin in the boundaryless case \cite{G}, the Birkhoff
normal form at a closed orbit is a spectral invariant. In
\cite{Z3} the author generalized this to the boundary case. Colin
de Verdiere's observation that the Birkhoff normal form at a
bouncing ball orbit determines an analytic plane domain with two
symmetries then gave the solution of the inverse spectral problem
for analytic plane domains with these symmetries \cite{Z3}. We
found it difficult however to calculate the normal form
coefficients by this method without assuming two symmetries. That
motivated us to try the Balian-Bloch approach.

 In the
recent paper \cite{ISZ}, Iantchenko-Sjostrand-Zworski give another
proof of the inverse  result with two symmetries  using   Birkhoff
normal forms of the billiard map and quantum monodromy operator
rather than the Laplacian; the method is quite elegant and
flexible. However, we believe that their normal form of the
monodromy operator is the same as the normal form in \cite{Z3} of
what is called there  the semiclassical wave group. We changed to
the present approach because we found it too difficult to
calculate the latter without the two symmetry assumption. However,
 the possibility remains that the Birkhoff normal form could be a better
 way to organize the inverse spectral data
than the wave invariants themselves.  The question is, how  much
of the boundary defining function is determined by the Birkhoff
normal form of the wave group or the monodromy operator at a
periodic orbit? It is known that the classical Birkhoff normal
form by itself will not determine the Taylor coefficients of the
domain unless there are two symmetries.

{\bf Acknowledgements.} We are heavily indebted to M.
Zworski for many helpful conversations on the contents of this
article, particularly during a visit at the Erwin Schrodinger
Institute during June 2001. In addition, we thank the ESI for
making the visit possible. We also thank J. Ralston for confirming
that various results of \cite{BGR} which we need are also valid in
the even-dimensional case. Finally, we thank the IMA for the
opportunity to present the results in their final form.

\section{Billiards and the length functional}\label{BACKGROUND}

Throughout this paper, we  will assume  that our obstacle is
up-down symmetric in the sense that there is an isometric
involution $\sigma$ of $\Omega$ which interchanges the endpoints
of a bouncing ball orbit (extremal diameter). However,
 we  only use the
assumption in the last step of calculating wave invariants and
determining the domain. Elsewhere, we make it for notational
convenience. We now extend the notation and terminology of
\cite{Z1} so that it applies to exterior domains which are
complements of obstacles with two  components.

We align the obstacles so that the bouncing ball orbit $\gamma$
projects to the line segment  $\overline{a_- a_+}$ along the
$y$-axis and so that its midpoint lies at the orgin  $0 \in \Rr^2$
and its  endpoints are   $a_{\pm} = (0, \pm L/2).$  We refer to
the top component of $\Omega$ as the $+$ component ${\mathcal
O}_+$ and the bottom one as the $-$ component ${\mathcal O}_-$.

There are two convenient parametrizations of $\partial \Omega.$
First, we denote by $q_{\pm } (\theta)$  the arc-length
parameterization of ${\mathcal O}_{\pm}$ with $q_{\pm}(0) =
a_{\pm}.$ At the end we will also use graph parametrizations:
 In a small strip $T_{\epsilon}(\gamma)$ around $\overline{a_1 a_2}$, the boundary consists
of two components which are symmetric graphs over the $x$-axis. We
write the graphs in the form $ y = \pm f(x)$ near $a_{\pm}$.

In the multiple reflection expansion of the exterior resolvent
$\RC(k + i \tau)$, we will encounter multiple integrals
$\int_{(\partial \Omega)^M}$. Since $\partial \Omega$ consists of
two components, each such term breaks up into $2^M$ multiple
integrals over the circle ${\bf T}$.  The terms can be enumerated
by maps $\sigma: \{1, \dots, M\} \to \{\pm \}$. We needed the same
enumeration in \cite{Z2} to denote the two local components of
$\Omega$ at the  ends of the bouncing ball orbit.

\subsection{Length functional}

For each map $\sigma : \{1, \dots, M\} \to \{\pm \}$ as above, we
define a length functional on ${\bf T}^M$ by:
\setcounter{equation}{3}
\begin{equation}\label{LENGTH}
\ba{lll}
\dis
L_{\sigma}(\phi_1, \dots, \phi_M) \\[6pt]
\dis \hspace{28pt}
=  |q_{\sigma_1} (\phi_1) - q_{\sigma_2} (\phi_2)| + \dots + |q_{\sigma_{M-1}} (\phi_{M - 1} ) -q_{\sigma_M} (\phi_M)| .
\ea
\end{equation}

It is clear  that $L_{\sigma}$ is a smooth function if $\sigma(j) \not= \sigma(j + 1)$ for any $j$. When
there do exist such $j$, then $L_{\sigma}$ is singular
on  the  `large diagonal' $\Delta_{j, j + 1}:= \{\phi_j = \phi_{j + 1}\}$, where it   has  $|x|$   singularities .  A standard (and easy) calculation shows that
\begin{equation}\label{Oneder} \begin{array}{rl}
\dis \frac{\partial}{\partial \phi_j} L _{\sigma} \ =\!\!&\dis
 \sin \angle ( q_{\sigma_{j + 1}}(\phi_{j + 1}) - q_{\sigma_j}(\phi_{j }),
\nu_{q(\phi_{j + 1})})\\[6pt]
&\dis -\ \sin \angle ( q_{\sigma_j}(\phi_{j }) - q_{\sigma_{j - 1}}(\phi_{j -1 }), \nu_{q_{\sigma_j}(\phi_j)}) .\end{array}\end{equation}
Here we  denote the acute angle between the link $q_{\sigma_{j + 1}} (\phi_{j + 1 }) - q_{\sigma_j} (\phi_{j })$ and
the inward unit normal $\nu_{q_{\sigma_{j + 1}}(\phi_{j + 1})}$ by $\angle ( q_{\sigma_{j +1}}(\phi_{j + 1}) - q_{\sigma_j}(\phi_{j }), \nu_{q_{\sigma_j}(\phi_j)}).$
 The condition that $\frac{\partial}{\partial \phi_j} L
= 0$ is thus that  the $2$-link defined by the triplet
$(q_{\sigma_{j-1}}(\phi_{ j- 1}, q_{\sigma_j}(\phi_j),
q_{\sigma_{j + 1}}(\phi_{j + 1})$ is  {\it Snell} at $\phi_j$,
i.e. satisfies the law of equal angles at this point. A smooth
critical point of $L$ on ${\bf T}^M$ is thus the same as an
$M$-link Snell polygon.   Note that Snell polygons include
polygons which have links inside the obstacles. They are known as
`ghost orbits' in the physics literature, since they cross the
boundary between  the interior and exterior domains. One may also
include singular critical points of $L_{\sigma}$ which correspond
to polygons in which an edge is collapsed.

\subsection{Length spectrum}

We will denote by $Lsp(\Omega^c)$, resp. $Lsp(\Omega)$,  the
lengths of periodic billiard trajectories in $\Omega^c,$ resp.
$\Omega.$  We do not include ghost orbits.  By periodic
trajectory, we mean closed bicharacteristic curve; for the
definition see \cite{GM, PS}. Special periodic trajectories are
$m$-link periodic reflecting ray,  i.e.  $m$-link Snell polygons.
Their lengths are isolated in the length spectrum.

When $\Omega^c$ is the complement of a pair of disjoint, strictly
convex obstacles, the only exterior transversal reflecting rays
are iterates $\gamma^r$ of the bouncing ball orbit $\gamma$
between the components, i.e. the line segment between the points
$a_{\pm} \in {\mathcal O}_{\pm}$. When ${\mathcal O}$ is not
assumed convex, there may be additional periodic reflecting rays
and the shortest periodic orbit of $\Omega^c$ need not be the
invariant bouncing ball orbit. Recall that we assume throughout
that $L_{\gamma}$ is isolated in $Lsp(\Omega) \cup Lsp(\Omega^c)$.

It is part of the folklore of the Balian-Bloch approach that
contributions to the trace formula due to ghost orbits are exactly
cancelled by other ghost orbits. We will not prove that this
cancellation occurs since the rigorous Poisson relation already
implies that only $\gamma^r$ contribute to the wave trace.
Additionally, there are so-called  diffraction orbits, which have
points of tangential intersection with the boundary and may glide
along the boundary.  They are not critical points of the length
functionals, and correspond not  to smooth singularities in the
wave trace, but to analytic singularities.

\section{Poisson relation for the exterior problem}
\label{POISSON}

We now recall the Poisson relation for
exterior domains in dimension 2.

The exterior {\it  Green's kernel} $G_B^{\Omega^c}(k + i
\tau, x, y) \in \mathcal D'(\Omega^c \times \Omega^c)$ with
boundary condition $B = D,$ or $ N$ (Dirichlet or Neumann) is the
kernel of the exterior  resolvent $\RB(k + i \tau) = -
(\DB +( k + i \tau)^2)^{-1} : H^{s}(\Omega^c) \to H^{s +
2}(\Omega^c)$. Here, $k \in \R$  and $\tau > 0$. Its kernel (the
exterior Green's function) may be characterized as  the unique
solution  of   the boundary problem:
\setcounter{equation}{5}
\begin{equation}\label{GREENK}
\left\{ \begin{array}{l}
\dis
- (\DB + (k + i \tau)^2) G_{B}^{\Omega^c}(k + i \tau, x, y) = \delta_{y}(x),\;\;\;
(x, y \in \Omega) \\ \\
\dis
B G_{B}^{\Omega^c}(k + i \tau, x, y) = 0, \;\;\; x \in \partial \Omega
\\ \\
\dis
 \frac{\partial G_{B}^{\Omega^c}(k \!+\! i \tau, x, y)}{\partial r} - i (k
\!+\! i \tau) G_{B}^{\Omega^c}(k \!+\! i \tau, x, y) = o(\frac{1}{r}) ,\;\; \mbox{as}\;\; r \to \infty . \end{array} \right.
\hspace{-15pt}
\end{equation}
Here, the boundary operator could be  either $B u = u|_{\partial
\Omega^c}$ ($D$) or $Bu =
\partial_{\nu} u |_{\partial
\Omega^c}$ ($N$). We use a similar notation for the interior
resolvent and boundary conditions. The boundary conditions $D$ and
$N$  are in a certain sense complementary, and we will write $B'$
for the complementary boundary condition to $B$.

 We are interested in the regularized
distribution trace of the combined operator $\RB(k + i \tau)
\oplus R^{\Omega}_{B'}(k + i \tau)$ , more usually regarded as a
distribution trace of the combined wave group $E_{B}^{\Omega^c}(t)
\oplus E_{B'}^{\Omega}(t)$, where $E_{B}^{\Omega^c}(t) = \cos t
\sqrt{- \DB})$. Let $\hat{\rho} \in C_0^{\infty}(\Rr^+)$ and
define
\begin{equation}
\ba{lll}\dis
R_{\rho B}^{\Omega^c}(k + i
\tau) \oplus R_{\rho B'}^{\Omega}(k + i \tau) \\[6pt]
\dis \hspace{32pt} = \int_{\R} \rho(k -
\mu) (\mu + i \tau) \RB(\mu + i \tau)\oplus R_{\rho
B'}^{\Omega}(\mu + i \tau) d\mu.
\ea
\end{equation}
From the resolvent
identity (e.g.)
$$\RB(\mu + i \tau) = \frac{1}{\mu + i \tau} \int_0^{\infty} e^{i (\mu + i \tau) t} E_{B}^{\Omega^c}(t) dt,\;\;\;\; $$
it follows that
\begin{equation} \label{RW} R_{\rho B}^{\Omega^c}(k + i \tau)  \oplus R_{\rho B'}^{\Omega}(k + i \tau)
= \int_0^{\infty} \hat{\rho}(t) e^{i (k + i \tau) t}
E_{B}^{\Omega^c}(t)\oplus E_{B'}^{\Omega}(t) dt. \end{equation}

\subsection{Poisson relation and scattering phase}

The following is the Birman-Krein formula for the relative trace
of the wave group of an exterior problem, as generalized by
Lax-Phillips and  Bardos-Guillot-Ralston (see \cite{BGR}, Theoreme
3) : Let $\hat{\phi} \in C_0^{\infty}(\R)$. Then:
$$Tr \int_{\R} \hat{\phi}(t) [E_B^{\Omega^c}(t)\oplus
E_{B'}^{\Omega}(t) - E_0(t)] dt = \frac{1}{2 \pi } \int_{\R}
\phi(\lambda) [d s_B(\lambda) + dN_{B'}(\lambda)]. $$ Putting
$\hat{\phi}(t) =  \hat{\rho}(t) e^{i (k + i \tau) t}$ and
rewriting in terms of the resolvent, we get: \begin{equation}
\label{BK}  \ba{lll} \dis Tr [R_{B \rho}^{\Omega^c} (k + i \tau)
\oplus R_{B' \rho}^{\Omega}
(k + i \tau) - R_{0 \rho}(k + i \tau)] \\[6pt]
\dis \hspace{35pt} =
 \frac{1}{2 \pi } \int_{\R} \rho(k + i \tau - \lambda)
[d s_B(\lambda)  + dN_{B'}(\lambda)].
\ea
\end{equation}
Here, $s_B(\lambda) = \log \det S_B(\lambda)$ where $S_B(E)$ is
the scattering operator and $N_{B'}$ is the interior Weyl counting
function with the indicated boundary condtions.

The singular support of the regularized  wave trace $Tr
[E_B^{\Omega^c}(t)  - E_0(t)] $  is contained in the set of
lengths of exterior periodic  generalized billiard trajectories.
That is,  by \cite{BGR}, Theoreme 5, we have
$$singsupp Tr [E_{B}^{\Omega^c}(t) - E_0(t)] \subset Lsp(\Omega^c).$$
Similarly for the interior wave group with $Lsp(\Omega)$ replacing
$Lsp(\Omega^c)$.
 When $L_{\gamma}$ is the length of a
non-degenerate periodic reflecting ray $\gamma$, and when
$L_{\gamma}$ is not the length of any other generalized periodic
orbit, then $Tr [E_{B}^{\Omega^c}(t)- E_0(t)]$ is a Lagrangian
distribution in the interval  $(L_{\gamma} - \epsilon, L_{\gamma}
+ \epsilon)$ for sufficiently small $\epsilon$,  hence $Tr
[E_{\rho B}^{\Omega^c}(k + i \tau) - E_{\rho 0}(t)]$ has a
complete asymptotic expansion in powers of $k^{-1}.$ Let us recall
the  statement in the exterior case  (see \cite{GM}, Theorem 1,
and  \cite{PS} Theorem 6.3.1 for the interior case, and
\cite{BGR}, \S 6 for the exterior case). Let $\gamma$ be a
non-degenerate billiard trajectory whose length $L_{\gamma}$ is
isolated and of multiplicity one in $Lsp (\Omega)$.  Let
$\Gamma_L$ be a sufficiently small conic neighborhood of $\Rr^+
\gamma$ and let $\chi$ be a microlocal cutoff to $\Gamma_L$.  Then
for $t$ near $L_{\gamma}$, the trace of the wave group has the
singularity expansion
\begin{equation} \label{TR}
\begin{array}{lll}\dis
\dis Tr \chi [ E_B^{\Omega^c} (t) \!-\! E_0(t)]\!\! & \sim\!\! &
a_{\gamma} (t \!-\! L_{\gamma} \!+\! i0)^{-1} \!+\! a_{\gamma 0}
\log (t \!-\! L_{\gamma}
\!+\! i 0) \\[4pt] & & \dis+
\sum_{k = 1}^{\infty}a_{\gamma k}(t \!-\! L_{\gamma} \!+\! i 0)^k \log (t
\!-\! L_{\gamma} \!+\! i 0) \end{array} \end{equation} where the coefficients $a_{\gamma k}$ are
calculated by the stationary phase method from a Lagrangian
parametrix. In the interior case, one of course omits the term $E_0(t)$.

The result may be re-stated as follows: Let $\hat{\rho} \in
C_0^{\infty}(L_{\gamma} - \epsilon, L_{\gamma} + \epsilon)$, equal
to one on $(L_{\gamma} - \epsilon/2, L_{\gamma} + \epsilon/2)$ and
with no other lengths in its support.  Then the interior trace
$Tr R_{B \rho}^{\Omega}(k + i \tau)$ and the exterior trace $Tr [R_{B\rho}^{\Omega^c}(k + i
\tau) - R_{0 \rho}(k + i \tau)]$ admit  complete asymptotic expansions of the form
\begin{equation} \label{ASYMP} \left\{ \begin{array}{l}
\dis Tr [R_{B\rho}^{\Omega^c}(k + i
\tau) - R_{0 \rho}(k + i \tau)] \sim \sum_{j = 0}^{\infty} B_{\gamma; j}
k^{-j} \\[9pt]
\dis Tr R_{B \rho}^{\Omega}(k + i
\tau) \sim e^{i (k + i \tau) \tau L_{\gamma}} \sum_{j = 0}^{\infty} B_{\gamma; j}
k^{-j}, \end{array} \right. \end{equation}
 whose  coefficients
$B_{\gamma; j}$ are canonically related to the wave invariants
$a_{\gamma; j}$ of periodic (internal, resp. external) billiard
orbits. We have removed the cutoff operator $\chi$ since  there
are no singularities of the trace at $t = L_{\gamma}$ in the
microsupport of $I - \chi$.  The coefficients depend on the choice
of boundary condition but we do not indicate  this in the
notation.

\subsection{Poisson relation and resonances}

We now recall the Poisson relation in dimension 2. Let
$\{\lambda_j\}$ denote the resonances of $\Delta_B^{\Omega^c}$,
i.e. the poles of the analytic continuation of $R_B^{\Omega^c} (k
+ i \tau)$ from $\{\tau > 0\}$ to the logarithmic plane $\Lambda.$
Let $\theta$ denote a (small) angle, and let $\Lambda_{\theta}$
denote the conic neighborhood of the real axis in $\Lambda$
consisting of points with angular coordinates in $(- \theta,
\theta).$ Further, let $\hat{\phi} \in C_0^{\infty}(\Rr^+)$ and
let $E_0(t)$ denote the free even wave kernel.  The global Poisson
formula (asserts the following:
\begin{prop} \label{ZWORSKI}  (see
Zworski \cite{Zw}, Theorem 1)
$$\begin{array}{c}
\dis Tr (\int_{\R} \hat{\phi}(t) (E_{B}^{\Omega^c} (t) - E_0(t))
dt = \sum_{\lambda_j \in \Lambda_{\theta}} m(\lambda_j)
\phi(\lambda_j)
+ m(0) \phi(0) \\[9pt] \dis
\hspace{110pt} +\ 2 \int_0^{\infty}\!\! \psi(\lambda)
\phi(\lambda) \frac{d s_B}{d \lambda} d \lambda +
\int_0^{\infty}\!\! \hat{\phi} (t)\nu_{\theta, \psi}(t) dt,\\[12pt]
\dis\mbox{with}\;\;\; \nu_{\theta, \psi} \in C^{\infty}(\R
\backslash \{0\}), \;\;\; \partial_t^{k} \nu_{\theta, \psi} =
O(t^{-N}),\;\; |t| \to \infty, \end{array} $$ where $s_B$
denotes the scattering phase of $\Delta_B^{\Omega^c},$ and where
$m(\lambda)$ denotes the multiplicity of the resonance $\lambda.$
Also, $\psi \in C_0^{\infty}(\R)$ is a cutoff which equals $1$ for
$t$ near $0$. \end{prop}

(We note that the sum over eigenvalues term in the formula stated
in \cite{Zw} ( Theorem 1) is absent   in the case of an exterior
domain.)

Substituting   $\hat{\phi}(t) = \hat{\rho}(t) e^{i(k + i \tau)
t}$, with $\hat{\rho} \in C_0^{\infty}(\Rr^+)$ with
supp$\hat{\rho}$ sufficiently close to $r L$ for some $r$,  into
(\ref{ZWORSKI})  and using (\ref{RW}), we obtain:
\begin{equation}\label{RESPOI}
\begin{array}{l}
\dis Tr [R_{\rho B}^{\Omega^c}(k + i \tau) - R_{0 \rho}(k + i \tau)]
\\[6pt]
\dis\hspace{50pt}
= \sum_{\lambda_j \in \Lambda_{\theta}} m(\lambda_j) \rho(\lambda + k + i \tau)
+  O(|k|^{-M}),\;\; \;\; (k \to \infty). \end{array} \end{equation}
In the remainder estimate we use that
$$
m(0) \rho(k \!+\! i \tau),\;\;\; \int_0^{\infty}\!\! \psi(\lambda)
\rho(\lambda \!+\! k \!+\! i \tau) \frac{d s_B}{d \lambda} d
\lambda,\;\;\;  \langle \nu_{\theta, \psi}, \hat{\rho}(t) e^{i(k \!+\! i \tau) t} \rangle
\!=\! O(k^{-\infty}). $$ In the case of the first two terms, this
follows immediately from the fact that $\phi(\lambda) = \rho(k + i
\tau - \lambda)$ and that both terms are integrals over a compact
set of values of $\lambda.$ For the third term we additionally use
that $\hat{\rho}$ vanishes in a neighborhood of $t = 0$;
integration by parts and the estimate on on $\nu_{\theta, \psi}$
then gives the rapid decay in $k$ of the term. It follows that the
asymptotics of the regularized resolvent trace $Tr [R_{\rho
B}^{\Omega^c}(k + i \tau) - R_{0 \rho}(k + i \tau)]$ are an
invariant of the resonance poles.

The following well-known proposition, explained to the author by   M. Zworski, implies more:  only poles in a logarithmic neighborhood  $\{|\Im \lambda_j| > N  \log |\lambda_j|\}$ of the
real axis cause singularities in the wave trace:

\begin{prop} For any $N \geq 0$, there exists $k(N)$ with $k(N) \to \infty$ as
$N \to \infty$, such that
$$\sum_{\lambda_j: |\Im \lambda_j| < \rho |\lambda_j|, |\Im \lambda_j| > N  \log |\lambda_j|}
e^{i \lambda_j t} \in C^{k(N)}( \R - \{0\}).$$ \end{prop}

\begin{proof} We have
\begin{equation}
\begin{array}{lll}
\dis | \partial_t^k \sum_{\lambda_j: |\Im \lambda_j| < \rho |\lambda_j|, |\Im \lambda_j| > N  \log |\lambda_j|}
e^{i \lambda_j t} |\\[19pt]
\hspace{100pt}
\dis\leq  \sum_{\lambda_j: |\Im \lambda_j| < \rho |\lambda_j|, |\Im \lambda_j| > N  |\lambda_j|}
|\lambda_j|^k |\lambda_j|^{-N t} \\[19pt]
\hspace{100pt} \dis \leq  \int_{\Rr^+} r^{k - t N} dN(r) \leq C
\int_{\Rr^+} r^{n + k - 1 - tN} dr.
\end{array} \end{equation}

In the last line we used a polynomial bound on the number of resonance poles
(cf.  \cite{Zw} (2.1) for the estimate).
Thus, the portion of the wave trace coming from poles outside of sufficiently
large logarithmic neighborhoods of the real axis is as smooth as desired.
\hfill\end{proof}

\section{Neumann expansion and  Poisson relation}\label{INTEXT}

In the Balian-Bloch approach to the Poisson relation, one uses the
classical theory of layer potentials (\cite{TII} \S 7. 11)   to
express $G_{D}^{\Omega^c}(k + i \tau, x, y)$ in terms of the
`free' Green's function  $G_0(k + i \tau, x, y)$ of $\Rr^2$, i.e.
the kernel of the free resolvent $(\Delta_0 - k + i \tau^2)^{-1}$
of the Laplacian $\Delta_0$ on $\Rr^2$. A classical reference for
this approach is the paper of Pleijel \cite{P}.  This approach is
also used in \cite{M} to make a reduction to the boundary in
studying resonance poles and in the duality result of (\ref{EP}).
We review the necessary background and notation.

\subsection{Potential theory}

We recall that
 the free Green's function in dimension two  is given by:
$ G_0(k + i \tau, x, y) =
  H^{(1)}_0(k + i \tau |x - y|)$, where    $H^{(1)}_0(z)$ is the Hankel function of index $0$.
We recall
its small and large distance asymptotics  (\cite{TI}, Chapter 3, \S 6):
\setcounter{equation}{12}
\begin{equation}\label{LK}
H^{(1)}_{\nu} (r) =  \left\{ \begin{array}{l}
\dis
 - \frac{1}{2\pi} \ln ( r) \;\rm{as}\; r \to 0, \; \mbox{if}\; \nu = 0\\ \\
\dis
 - \frac{i \Gamma(\nu)}{\pi} (\frac{2}{r })^{\nu}  \;\rm{as}\; r \to 0,  \; \mbox{if}\; \nu > 0\\ \\
\dis
  \frac{e^{i (k + i \tau r - \nu \pi/2 - \pi/4)}}{r^{1/2}}
\bigg(1 +   O\Big(\frac{1}{|r|^3}\Big)\bigg) \;\rm{as}\;
  |k + i \tau| r \to \infty.  \end{array} \right.
\hspace{-10pt}
\end{equation}
Here it is assumed that $\tau> 0$. Thus,    $ G_0(k + i \tau, x,
y)$ has two kinds of asymptotics: a semiclassical asymptotics for
large $(k + i \tau) |x - y|$ and a singularity asymptotics along
the diagonal. It is the first kind of asymptotics which accounts
for the connection to billiard trajectories. The singularity along
the diagonal will have to be regularized.

\subsubsection{Layer potentials and boundary integral operators}

The double layer potential is   the
operator
\begin{equation} \label{layers}
{\mathcal D} \ell(k + i \tau)f(x) =  \int_{\partial
\Omega}\frac{\partial}{\partial \nu_y} G_0(k + i \tau, x, q) f(q)
ds(q),
 \end{equation}
from $H^s(\partial \Omega) \to H^{s+ 1/2}_{loc} (\Omega)$, where
$ds(q)$ is the arc-length measure on $\partial \Omega$, where
$\nu$ is the interior unit normal to $\Omega$, and  where
$\partial_{\nu} = \nu \cdot \nabla$. It  induces the boundary
operator
\begin{equation} \label{blayers} N(k + i \tau)f(q) =  2 \int_{\partial
\Omega}\frac{\partial}{\partial \nu_y} G_0(k + i \tau, q, q')
f(q') ds(q')
\end{equation} which map $H^s(\partial
\Omega) \to H^{s + 1}(\partial \Omega).$

By the explicit formula we have:
$$ \begin{array}{lll}
\dis \frac{1}{2} N(k + i \tau, q(\phi), q(\phi')) & = &
\dis \partial_{\nu_y} G_0(\mu,  q(\phi), q(\phi')) \\ [6pt]
& = &\dis
-\ (k + i \tau ) H^{(1)}_1 (k + i \tau|q(\phi) - q(\phi')|)\\ [6pt]
& &\times\ \dis \cos \angle(q(\phi) - q(\phi'),  \nu_{q(\phi)}).
\end{array}$$

We will need the `jump' formula of potential theory (
\cite{TI} Chapter 7. 11, (11.7)),
\begin{equation} \label{BV} ({\mathcal D}\ell(k + i \tau) u)_{\pm}(q) = \frac{1}{2} (N(k + i \tau) \mp I) u(q), \end{equation}
where
\begin{equation} \left\{ \begin{array}{l} f_+(x) :=  \lim_{x \to q, x \in {\mathcal O}} f(x), \\ \\
f_-(x) = \lim_{x \to q, x \in \Omega} f(x). \end{array} \right.
\end{equation}

Note that the sign of each term in (\ref{BV}) depends on the
choice of the interior/exterior: in $N(k + i \tau)$ it depends on
the choice of interior/exterior unit normal and in the identity
operator term it depends on whether the limit is taken from the
interior or exterior.  In the interior case we have $({\mathcal
D}\ell(k + i \tau) u)_{\pm}(q) = \frac{1}{2} (N(k + i \tau) \pm I)
u$ where $N(k + i \tau)$ is defined using the interior normal and
where $\pm$ have the same meanings as in the interior case.

\subsection{Interior/Exterior duality}

In \cite{Z1, Z2} we used the classical reduction of the Dirichlet
problem to the boundary  to study the Dirichlet resolvent (as in
\cite{BB1, BB2}). We now use a similar method to obtain a
convenient formula for the Fredholm determinant, $ \det ( (I +
N(k))).$ Similar ideas can be found in \cite{THS, THS2,  AG, GP}.

The following formula is sometimes referred to as interior/\break
exterior duality. Combined with the Birman-Krein formula (\ref
{BK}) it expresses the growth rate of the inside plus outside
spectra in terms of the determinant of a boundary integral
operator:

\begin{prop} \label{MAIN} For any $\tau \geq 0$, the  operator $ (I + N(k + i \tau))$ is of trace class
and has a well-defined Frehdolm determinant, and we have:
$$\begin{array}{l}
\dis Tr_{\Rr^2} [R_{ \rho D}^{\Omega^c}(k + i \tau) \oplus
R_{N \rho}^{\Omega}  - R_{0 \rho}(k + i \tau)]\\[6pt]
\dis\hspace{43pt} = \int_{\R} \rho(k - \lambda)  \frac{d}{d
\lambda}  \log \det (I + N(\lambda + i \tau)
d\lambda,.\end{array}$$ where $\det (I + N(\lambda + i \tau)$ is
the Fredholm determinant.

\end{prop}

\begin{proof}
We first argue formally. The interior/exterior resolvent kernels
can be constructed in the classical way, as follows (see e.g.
\cite{P}). For simplicity, we consider the interior Dirichlet
resolvent, but the construction for the exterior resolvent or with
a change to Neumann boundary conditions is almost the same. We
have:
\begin{equation} \label{NS} \begin{array}{l} R^{\Omega^c}_D(k \!+\! i \tau) =
1_{\Omega^c} [ R_0(k \!+\! i \tau)  \!-\! {\mathcal D} \ell(k
\!+\! i \tau) (I \!+\!  N(k
\!+\! i \tau))^{-1} \gamma R_0(k \!+\! i \tau)] 1_{\Omega^c}, \\ \\
(R^{\Omega}_N(k \!+\! i \tau))^{tr}  = 1_{\Omega}[
 R_0(k \!+\! i \tau)\!-\! {\mathcal D} \ell(k \!+\! i \tau) (I \!+\!  N(k
\!+\! i \tau))^{-1} \gamma R_0(k \!+\! i \tau)] 1_{\Omega} ,
\end{array}
\end{equation}
where $\gamma$ denotes the restriction to the boundary taken from
within the relevant domain.  When we take the regularized trace,
we first subtract from $R_N^{\Omega}(k + i \tau) \oplus
R_D^{\Omega^c}(k + i \tau)$ the free operator $R_0(k + i \tau)$,
which removes the first terms on the right.  We then cycle the
factor $\gamma R_0(k + i \tau)$ from the right to the left side,
obtaining an operator on $\partial \Omega.$ The inside and outside
terms add up to the kernel
\begin{equation} \label{KERNEL} \int_{\Rr^2} G_0(k + i \tau, q', y) \partial_{\nu_q} G_0(k+ i
\tau, y, q) dy. \end{equation} This indeed is why the interior
Dirichlet and exterior Neumann problems were combined and explains
the sense in which they are complementary.  We now claim that the
resulting kernel
 equals $\frac{1}{2 k} \frac{d}{d k} N(k + i \tau, q', q).$ To prove this we
note  that
$$\frac{d}{d \lambda} R_0(\lambda) = 2 \lambda (\Delta + \lambda^2)^{-2} = 2 \lambda
R_0(\lambda)^2.$$
Hence,
$$\frac{d}{d \lambda}N(\lambda) = \frac{d}{d\lambda}  \gamma \partial_{\nu_q}
R_0(\lambda) \gamma = 2 \lambda \gamma \partial_{\nu_q} R_0(\lambda) \circ R_0(\lambda) \gamma$$
which is precisely the kernel (\ref{KERNEL}). Hence the right side
of (\ref{NS}) equals $\frac{1}{k + i \tau} \frac{d}{d k } \log (I
+ N(k + i \tau)).$  The final formula follows by putting $k =
\lambda,$ integrating against $\rho(k - \lambda) (\lambda + i
\tau)$ and noting the cancellation of the second factor.

To justify the formal  manipulations, we need to show that the
relevant traces and  determinant are well defined. This will be
done in a sequence of Lemmas.

\begin{lem} For any $\tau$, $N(k + i \tau) \in {\mathcal I}_1(\partial \Omega)$, the ideal of trace
class operators on $L^2(\partial \Omega)$. Hence, $\det (I + N(k +
i \tau))$ is well-defined as a Fredholm determinant.
\end{lem}

\begin{proof}

In \cite{Z1} it is proved that  $N(k + i \tau) \in
\Psi^{-2}(\partial \Omega)$ and this immmediately implies that it
is of trace class. In fact, $N(k + i \tau, q, q')$ has just a $|q
- q'| \log |q - q'|$ singularity on the diagonal. That it is trace
class also follows from the Hille-Tamarkin theorem (see \cite{Z1}
for references).

It is then a classical remark that the Fredholm determinant $\det
(I + N(k + i \tau))$ is well-defined, in fact $|\det (I + N(k + i
\tau))| \leq \exp(||N(k + i \tau)||_1),$ where $|| \cdot ||_1$ is
the trace norm. See \cite{S}, Lemma 3.3.

\end{proof}

\begin{rem}

\noindent{\bf 1} Note that the usual statement (cf. \cite{TII}) is
that $N(k + i \tau)$  is of order $-1$. However, the principal
symbol vanishes in dimension $2$.

\noindent{\bf 2} We recall (see \cite{S}, Theorem 3.10) that the
Fredholm determinant is given by:
$$\det
(I + N(k + i \tau)) = \sum_{n=0}^{\infty} \frac{1}{n!}
\int_{{\partial \Omega}^n} \det[ N(k + i \tau, q_i, q_j)]_{1 \leq
i,j \leq n} d s(q_1) \cdots ds(q_n). $$
\end{rem}

Again combining classically known facts, we have:

\begin{lem} For any $\tau > 0, \log \det N(k + i \tau)$ is well-defined and  differentiable in $k$,  $ (I
+ N(k + i \tau))^{-1} N'(k + i \tau)$ is of trace class and we
have:
$$\frac{d}{d k} \log \det N(k + i \tau) = Tr_{\partial \Omega} (I
+ N(k + i \tau))^{-1} N'(k + i \tau). $$
\end{lem}

\begin{proof} First we note that $(I + N(k + i \tau))$ is
invertible on $L^2(\partial \Omega)$ if $\tau > 0$, so $\det (I +
N(k + i \tau)) \not= 0.$ Hence its logarithm is well-defined.
Differentiability of $\log \det (I + N(k + i \tau))$ and the
formula for the derivative follows from the general fact that
$\det(I + A)$ is Frechet differentiable on ${\mathcal I}_1$ with
derivative equal to $(I + A)^{-1}$ if $-1 \notin \sigma(A)$ (here,
$\sigma(A)$ is its spectrum; see \cite{S}, Corollary 5.2).

The singularity on the diagonal of $N(k + i \tau, q,q')$ is
independent of $k + i \tau$, so the derivative is at least as
regular. Hence, $N'(k + i \tau) \in {\mathcal I}_1.$ The statement
that $(I + N(k + i \tau))^{-1} N'(k + i \tau) \in {\mathcal I}_1$
follows from the fact that  $(I + N(k + i \tau)^{-1} \in
\Psi^0(\partial \Omega)$.

\end{proof}

\begin{lem} The operators
$$ \left\{ \begin{array}{ll} (i) & 1_{\Omega} \int_{\R} \rho(k - \lambda) {\mathcal D} \ell(\lambda \!+\! i \tau) (I
\!+\!  N(\lambda \!+\! i \tau))^{-1} \gamma R_0(\lambda \!+\! i
\tau)] d \lambda 1_{\Omega}, \\ & \\(ii) &
 1_{\Omega^c} \int_{\R} \rho(k - \lambda) {\mathcal D} \ell(\lambda \!+\! i \tau) (I \!+\!
 N(\lambda
\!+\! i \tau))^{-1} \gamma R_0(\lambda \!+\! i \tau)] d\lambda
1_{\Omega^c} ,
\end{array}\right.$$
are of trace class and the sum of their traces equals $Tr_{\Omega}
(I + N(k + i \tau))^{-1} \; N'(k + i \tau). $
\end{lem}

\begin{proof}

We expand  $(I \!+\!  N(\lambda \!+\! i \tau))^{-1} = I -
N(\lambda) (I \!+\!  N(\lambda \!+\! i \tau))^{-1}. $ The $I$ term
gives us ${\mathcal D} \ell(\lambda \!+\! i \tau)  \gamma
R_0(\lambda \!+\! i \tau)]$. We break up each of ${\mathcal D}
\ell(\lambda \!+\! i \tau)$ and  $\gamma R_0(\lambda \!+\! i
\tau)]$ into regular (continuous) and singular parts using the
small $z$ expansion (\ref{LK})  of the Hankel functions $H^1_0(z),
H^1_z(z)$. The singular parts of the kernels are  independent of
$(k + i \tau)$.  Hence the product of the singular parts  gets
multiplied by $\hat{\rho}(0) = 0$ upon integration against $\rho(k
- \lambda)$. Therefore only the regular parts survive the
convolution, and these are  of trace class. The remaining term is
a bounded operator composed with trace class operator for each
$\lambda$, so it is of trace class for each $\lambda$.

We may then cycle the factor of ${\mathcal D} \ell(\lambda \!+\! i
\tau) $ (or each of its regular and singular parts) to the right
side of the traces. After doing so, we  reassemble the regular and
singular parts into $   \gamma R_0(\lambda \!+\! i \tau)]{\mathcal
D} \ell(\lambda \!+\! i \tau) $and we reassemble $I - N(\lambda)
(I \!+\! N(\lambda \!+\! i \tau))^{-1}$ into $(I \!+\! N(\lambda
\!+\! i \tau))^{-1}$. The calculation at the beginning of the
proof then gives the stated formula.

\end{proof}

This concludes the proof of Proposition \ref{MAIN}.

\end{proof}

\begin{cor} \label{MAINCOR} Suppose that $L_{\gamma}$ is the only length
in the support of $\hat{\rho}$. Then,
$$\int_{\R} \rho(k -
\lambda) \frac{d}{d \lambda} \log \det (I + N(\lambda + i \tau)) d
\lambda  \sim \sum_{j = 0}^{\infty} B_{\gamma; j} k^{-j}, $$ where
$B_{\gamma; j}$ are the wave invariants of $\gamma$ in
(\ref{ASYMP}).

 \end{cor}

 To  prove Theorems (\ref{ONESYM})
and (\ref{ONESYMINT}), it thus suffices to determine $\Omega$ from
the integrals in Corollary (\ref{MAINCOR}).  Since we have
combined the interior/exterior, we emphasize that our inverse
result only assumes knowledge of the resonances: Indeed, suppose
the resonance poles of $\Omega^c$ are known. Then by Proposition
(\ref{ZWORSKI}), the asymptotics of the exterior resolvent trace
in terms of exterior periodic orbits are known. But by Proposition
(\ref{MAIN}), these are the same as the asymptotics of the
integrals in  Corollary (\ref{MAINCOR}).

\section{Trace asymptotics}

We now explain how to use  Corollary (\ref{MAINCOR}) to  calculate
the coefficients $B_{\gamma, j}$. We write  $$\frac{d}{d \lambda}
\log \det (I + N(\lambda + i \tau)) = Tr_{\partial \Omega} (I +
N(\lambda + i \tau))^{-1} N'(k + i \tau),$$ and then
 expand  $(I +  N(k + i \tau))^{-1}$ in a finite geometric
series plus remainder:
\begin{equation} \label{GS} (I \!+\! N(\lambda \!+\! i \tau))^{-1} = \sum_{M = 0}^{M_0} (-1)^M \; N(\lambda)^M
+
 (-1)^{M_0 + 1} \; N(\lambda)^{M_0 + 1}  (I
\!+\! N(\lambda \!+\! i \tau))^{-1}. \end{equation} We now argue
that for each  order $k^{-J}$ in the trace expansion of Corollary
(\ref{MAINCOR}) there exists $M_0(J)$ such that
\begin{equation} \label{ASY} \begin{array}{ll} (i) &
\sum_{M = 0}^{M_0} (-1)^M  Tr  \int_{\R} \rho(k - \lambda)\;
N(\lambda)^M N'(k + i \tau) d \lambda \\ & \\ &  = \sum_{j =
0}^{J} B_{\gamma; j} k^{-j} + O(k^{-J - 1}) , \\&  \\(ii) & Tr
\int_{\R} \rho(k - \lambda)  N(\lambda)^{M_0 + 1}  (I \!+\!
N(\lambda \!+\! i \tau))^{-1} N'(k + i \tau) d \lambda = O(k^{-J -
1}).
\end{array} \end{equation}

 We will sketch the proof of (i) in some detail since it explains
 how to calculate the coefficients $B_{\gamma^r, j}$. For a
 discussion of the remainder (in an analogous but not identical
 calculation) we refer to \cite{Z2}.

 To simplify the notation, we integrate by parts in (i) to throw
 the derivative onto $\rho$, and then the issue is to analyse the
 traces:
\setcounter{equation}{22}
 \begin{equation}\label{MINT}
\begin{array}{lll}
\dis \hspace{-25pt} Tr  \int_{\R}  \rho(k - \lambda) N(\lambda + i \tau)^M d \lambda \\[9pt]
\dis\hspace{0pt}  = \int_{\R} \int_{(\partial \Omega)^{M}}\rho(k -
\lambda)  [\Pi_{j = 1}^M N(\lambda \!+\! i \tau, q_j, q_{j\!+\!1})
ds(q_j)] d\lambda
\;\;\;(\mbox{where}\; q_{M \!+\! 1}\! =\! q_1);  \\[12pt]
\dis\hspace{0pt}
\end{array}
\hspace{-25pt}
\end{equation}

Since $$N(k + i \tau): L^2({\mathcal O}_+) \oplus  L^2({\mathcal
O}_-) \to  L^2({\mathcal O}_+) \oplus  L^2({\mathcal O}_-)$$ we
write it as
$$N(k)   = \left(
\begin{array}{ll}
 N_{++}(k)  & N_{+ -}(k)\\ & \\
N_{-+}(k)  & N_{--} (k)\end{array} \right). $$.

Then (\ref{MINT}) is a sum of terms  $\sum_{\sigma: \{1, \dots, M
\} \to \{\pm 1\} } I_{M, \rho}^{\sigma}$ with
 \begin{equation}\label{MINT2}
 I_{M, \rho}^{\sigma} = \int_{\R}  \int_{{\bf T}^{M }} \rho(k
 - \lambda)
[\Pi_{j = 1}^M N_{\sigma_j, \sigma_{j+1}} (\lambda \!+\! i \tau,
q_{\sigma_j}(\phi_j), q_{\sigma_{j\!+\!1}}(\phi_{j \!+\! 1}) d
\phi_j] d \lambda. \end{equation} Here and hereafter, indices are
understand modulo $M$. Using the asymptotics of the free Green's
function in (\ref{LK})-(\ref{HSYMBOL}), the $\sigma$th term is
formally an oscillatory integral with phases given by the length
functional (\ref{LENGTH}). We would like to calculate the trace
asymptotically by the stationary phase method.

As discussed extensively in \cite{Z1},  we cannot immediately do
so because of the singularities along the diagonals of the
integrand. We therefore first need to de-singularize the
integrals.  One of the main advantages of the present
inside/outside reduction to the boundary is that it simplifies the
regularization procedure by eliminating the integral over $\Omega$
in \cite{Z1}.   We briefly outline the method and refer the reader
to \cite{Z1} for further details.

\subsection{Boundary integral operators as quantized billiard maps}

The operators $N_{- +}(k + i \tau), N_{+ -}(k + i \tau)$ are
semiclassical  Fourier integral operators  with phases $|q_+(\phi)
- q_-(\phi')|$. The phase (and amplitude) is non-singular since
the boundary components are disjoint. The phase in the $+ -$ case
parametrizes the graph of the following multi-valued billiard map
from $B^*{\mathcal O}_+$ to $B^*{\mathcal O}_-$: given $(q_+, v_+)
\in T_q {\mathcal O}_+$ with $|v| < 1$, add a multiple of the unit
outward normal to turn $v$ into an outward unit vector and proceed
along the straight line in that direction and let $q_-$ be a point
of intersection of this line with  ${\mathcal O}_-$.  Let $v_-$ be
the tangential projection of the terminal velocity vector at $q_-$
, and put $\beta_{+ -}(q_+, v_+) = (q_-, v_-).$ We have
(deliberately) described $\beta_{+ -}$ in an ambiguous way: when
$q_-$ is the first intersection point, we have the usual exterior
billiard map; but the phase actually parametrizes the canonical
relation which includes intersection points which occur after the
line enters the interior of ${\mathcal O}_-$. As mentioned above,
these are `ghost orbits' of the billiard flow and they cancel out
of the trace formula and do not contribute to the wave invariants
we are calculating.

The operators $N_{++} (k + i \tau)$ and $N_{--}(k + i\tau)$ are
more complicated, and also, as it turns out, less important. Each
is a combination of a homogeneous pseudodifferential operator of
order $-1$ with singularity on the diagonal and a Fourier integral
operator of order $0$ which `quantizes' the billiard map of the
interiors of  $\Omega$ and $\Omega^c$. They are exactly the kind
of operators discussed extensively in   \cite{Z1, HZ}. There is a
further transition region between these two regimes in which it
behaves like an Airy operator, but this region will not be
important for our problem.

To separate the two basic regions (tangential and transverse) of
$N_{++}$ and $N_{--}$, we introduce a cutoff of the form
$\chi(k^{1 - \delta}(\phi - \phi'))$ in terms of the arc-length
coordinate $\phi$ used in the parameterization $q(\phi)$. Here,
$\chi$ a smooth bump function which is supported on $[-1,1]$ and
equals $1$ on $[-1/2, 1/2]$. We then write (dropping the
subscripts)
\begin{equation}\left\{ \begin{array}{l} N(k + i \tau) = N_0(k + i \tau) + N_1(k + i
\tau),\\ \\  N_0 (k + i \tau, q(\phi), q(\phi')) = \chi(k^{1 -
\delta}(\phi - \phi')) N (k + i \tau, q(\phi), q(\phi')),\\ \\ N_1
(k + i \tau, q(\phi), q(\phi')) = (1 - \chi(k^{1 - \delta}(\phi -
\phi')) N (k + i \tau, q(\phi), q(\phi')). \end{array} \right.
\end{equation}

Each term is a (non-standard) Fourier integral operator with
amplitudes belonging to  the following symbol class: We denote by
$S^p_{\delta}({\bf T}^m)$ the class of symbols $a(k, \phi_1,
\dots, \phi_m)$ satisfying:
\begin{equation}\label{SYMBOL}
|(k^{-1} D_{\phi})^{\alpha} a(k, \phi)| \leq C_{\alpha} |k|^{p -
\delta |\alpha|},\;\;\; (|k| \geq 1).  \end{equation} This follows
from the  asymptotics of Hankel functions (\ref{LK}), which  may
be described as follows:   there exist amplitudes $a_0, a_1$
satisfying
$$
\ba{ll} \dis (1 - \chi(k^{1 - \delta} z)) a_1((k + i \tau)z) \in
S^{0}_{\delta}
(\R),\\[6pt]\dis
   (1 - \chi(k^{1 - \delta} z))  a_0((k + i \tau)z)   \in
S^{-1/2}_{\delta} (\R) \ea $$ and such that

\begin{equation} \label{HSYMBOL} \left\{\!\!\!\begin{array}{ll}
 (i) &\dis
\!\!  H^{(1)}_0( (k + i \tau) z) = e^{i((k + i \tau) z } a_0((k + i \tau)z), \\ & \\
(ii) &\dis \!\!  (k \!+\! i \tau) H^{(1)}_1( (k \!+\! i \tau) z) =
(k \!+\! i \tau)^{ \frac{1}{2}} e^{- i((k \!+\! i \tau) z} a_1((k
\!+\! i \tau)z), \end{array}\right. \hspace{-15pt}\end{equation}

The $N_0$ term has a singularity on the diagonal of a
pseudodifferntial operator of order $-2$. The $N_1$ term is
manifestly an oscillatory integral operator of order $0$  with
phase $|q(\phi) - q(\phi')|$. As discussed in the $-+$ case, it
generates the interior billiard map (interpreted in the
multivalued sense of the $-+$ case).

\subsection{Microlocalization to $\gamma$}

So far, our considerations have been global. But our goal is to
calculate wave invariants at a bouncing ball orbit, and for this
purpose we can microlocalize the trace  to a neighborhood of
$\gamma$. This obvious sounding statement is somewhat problematic
in the present approach since the relevant operators are not
standard pseudodifferential or Fourier integral operators.

The billiard map $\beta_{+ -}$ is a cross section of the billiard
flow, and in this cross section $\gamma$ corresponds to the
periodic points $(0, 0) \in B^* {\mathcal O}_+$ and $(0,0) \in
B^*({\mathcal O}_-)$  of period $2$ for $\beta_{+ -}$. Here, we
choose the parametrizations of ${\mathcal O}_{\pm}$ such that
$\phi = 0$ at each endpoint of the bouncing ball orbit; it is
normal to ${\mathcal O}_{\pm}$, so its tangential part is the
vector $(0,0)$. To microlocalize to this orbit we introduce a
(block diagonal) semiclassical pseudodifferential cutoff operator
$$\chi(\phi, k^{-1} D_{\phi})  = \left(
\begin{array}{ll}
 \chi_{++} & 0\\ & \\
0 & \chi_{--} \end{array} \right). $$ with complete symbol
$\chi(\phi, \eta)$ supported in $V_{\epsilon}:= \{(\phi, \eta):
|\phi|, |\eta| \leq \epsilon\}$. We also write
$$ \rho * \frac{d}{d\lambda} \log
\det (I + N(k + i \tau )) = Tr \rho *  (I + N(k + i \tau ))^{-1}
\circ \frac{d}{dk} N(k + i \tau). $$

 Recall that to test whether $(L, \tau)  \in WF(u)$ (wave front set) for $u(t) \in {\mathcal
 S}'(\R)$, one checks whether ${\mathcal F} (\hat{\rho } u)(k) \sim 0$
 for all
 $\hat{\rho} \in C_0^{\infty}(\R)$ satisfying $\hat{\rho }\equiv 1
 $ in some interval $(L - \epsilon, L + \epsilon).$ Here,
 ${\mathcal F} u = \hat{u}$ denotes the Fourier transform.
 Equivalently, one checks whether $\rho * \hat{u} \sim 0$.
We now show  that in calculating wave invariants we can
microlocalize our non-standard operators to $\gamma$ in the
standard way.

As a preliminary step, we prove the simpler  relation
\begin{equation}\label{CUTOFF2}  \rho * Tr [{\mathcal D}(k + i \tau )
(I + N(k))^{-1 tr} {\mathcal S}(k + i \tau) ^{tr} \circ (1 -
 \tilde{\chi}_{\gamma}(k)] \sim 0. \end{equation}
 where $\tilde{\chi}_{\gamma}(k)$ is a semiclassical cutoff in
$\Omega^c$ to a microlocal neighborhood of $\gamma$ which has the
form
$$\tilde{\chi}_{\gamma}(k) = \tilde{\chi}(r, y, k^{-1} D_y) $$
near the boundary.  Here,  $(r, y)$ are Fermi normal coordinates
in $\Omega^c$ near the endponts of $\gamma$ , with $r$ the
distance to the boundary. Thus, differentiations are only in
tangential directions. Note that (\ref{CUTOFF2}) is equivalent to
saying that
 \begin{equation}\label{CUTOFF}  Tr[ R_{\rho D}^{\Omega^c} - 1_{\Omega^c} R_{0 \rho}] \circ (I -
 \tilde{\chi}_{\gamma}(k) )\sim 0, \end{equation}

 In view of  (\ref{NS}), (\ref{CUTOFF})-(\ref{CUTOFF2}) follow from standard
 wave front set ($WF$)  considerations for  the exterior even
wave group $E_D^{\Omega^c}(t, x, y)$, and from our assumption that
$\gamma^r$ is the only exterior periodic orbit of its length.
Recall that $WF(E_D^{\Omega^c}(t, x, y))$ is the space-time graph
\begin{equation} \Gamma =  \{(t, \tau, x, \xi, x', \xi') \in T^*(\R \times \Omega^c \times
\Omega^c):\;  \tau = - |\xi|, \;\; G^t(x, \xi) = (x', \xi')\},
\end{equation}
of the generalized billiard flow in the exterior. The only unusual
feature is our particular choice of cutoff operator to
microlocalize to $\gamma$. To verify that it behaves correctly
like a microlocal cutoff to $\gamma$ we write
\begin{equation}  R_{\rho D}^{\Omega^c}(k )
\tilde{\chi}_{\gamma}(k) = \int_0^{\infty} \hat{\rho}(t)
E_{D}^{\Omega^c}(t) \tilde{\chi}_{\gamma}(r, y, |D_t|^{-1} D_y)
 e^{i (k + i \tau ) t} dt. \end{equation}
 Now in Fermi coordinates $(r, \rho dr, y, \eta dy), \Rr^+ \gamma$ is ray in the direction
 of $dr$. Since the space-time
 graph of  the cotangent bundle along $\gamma$ may be described in normal coordinates as a
 neighborhood of
$$\{(t, \tau, t, \tau, 0, 0)  \}$$
a conic neighborhood  may be described in these coordinates  by
$|y| \leq \epsilon, |\eta/\tau| \leq \epsilon$. This is precisely
the set to which $\tilde{\chi}_{\gamma}(r, y, |D_t|^{-1} D_y)$
microlocalizes. Thus, (\ref{CUTOFF}) and hence (\ref{CUTOFF2}) are
 correct.

We now show that  this  empty WF relation remains true after
moving the cutoff to the left of $\SL^{tr}$ as a cutoff
pseudodifferential operator on the boundary to a neighborhood of
the billiard map orbit. As mentioned above, this would be
straightforward if $\DL$ and the other operators in the
composition were  standard FIO's. However in addition to their
semiclassical canonical relations they have singularities on the
diagonal and both contribute to their $WF$. We will shorten the
notation for the cutoff operator to $\chi_0$.

 \begin{lem} We have:
  $$\begin{array}{l}Tr \rho *  (I + N(k + i \tau ))^{-1} \circ \frac{d}{dk} N(k
+ i \tau) \\ \\ \sim Tr \rho *  (I + N(k + i \tau ))^{-1} \circ
\chi_0(k) \circ \frac{d}{dk} N(k + i \tau) .\end{array} $$

 \end{lem}

 \begin{proof}
We first express the various layer potentials in terms of the wave
group:
$$\left\{ \begin{array}{l} {\mathcal S}\ell(k + i \tau )^{tr} = \gamma  G_0(k + i \tau )
= \int_0^{\infty} e^{i t (k + i \tau ) }  \gamma E_0(t) dt,\\ \\
 1_{\Omega^c} ({\mathcal D}\ell(k + i \tau) \circ (I + N(k + i \tau ))^{-1 }(x, q)
= P_N^{\Omega^c}(k + i \tau ,x, q) \\ \\
= \gamma_{\nu 2} G_N^{\Omega^c}(k  + i \tau, x, q) =
  \int_0^{\infty} e^{i t (k + i \tau )}  \gamma_{\nu 2} E_{N}^{\Omega^c}(t)dt ,\end{array}
\right.$$ where $P_N^{\Omega^c}$ is the Neumann Poisson kernel of
the exterior domain and where  $\gamma_{\nu} u =
\partial_{\nu} u|_{\partial \Omega^c}, \gamma u = u|_{\partial \Omega^c}.  $   We have subscripted the restriction operators to clarify which
variables they operate on. The composition $ \DL \circ (I +
N(k))^{-1  }  \circ (1 - \chi_{\gamma}(k)) \circ \SL^{tr}$ may be
written (with the relevant value of $k$)  as
\begin{equation} \int_0^{\infty} e^{i t (k + i \tau )} \{\int_0^{t} \gamma_{\nu 2}
E_{N}^{\Omega^c} (t - s) \circ (I - \chi_{\gamma})(y, |D_t|^{-1}
D_y)
 \gamma  \circ  E_0(s) ds\} dt. \end{equation}

 We therefore have:
\begin{equation}\label{WFUB2}\begin{array}{l}  \rho *   \SL \circ (I +
N(k))^{-1 tr }  \circ (1 - \chi_{\gamma}(k)) \circ \DL^{tr}
\\ \\
= \int_0^{\infty} \hat{\rho}(t) [\int_0^t \gamma_{\nu 2}
E_{N}^{\Omega^c} (t - s) \circ (I - \chi_{\gamma})(y, |D_t|^{-1}
D_y)
 \gamma_{\nu 1} \circ  E_0(s) ds] e^{i (k + i \tau) t}
 dt .
\end{array}
\end{equation}

Our goal is to show that
\begin{equation}
WF [Tr \int_0^t\gamma_{\nu 2} E_{N}^{\Omega^c} (t - s) \circ (I -
\chi_{\gamma})(y, |D_t|^{-1} D_y)
 \gamma \circ  E_0(s) ds] \cap [rL_{\gamma} - \epsilon, r
 L_{\gamma} + \epsilon] = \emptyset. \end{equation}
This follows as long as  $$V(s, t) =   Tr \gamma_{\nu 2}
E_{N}^{\Omega^c} (t - s) \circ (I - \chi_{\gamma})(y, |D_t|^{-1}
D_y)
 \gamma \circ  E_0(s) $$
 is a smooth function for $t \in (r L_{\gamma} - \epsilon,
 rL_{\gamma} + \epsilon)$ and for $s \in (0,  rL_{\gamma} + \epsilon)$.
 To prove that it is smooth, we observe that the singular support consists of $(s, t)$  such that there
 exists a closed billiard orbit of length $t$ outside a phase space neighborhood
 of $\gamma^r$ which consists of a
 straight line segment of length $s$ from a point $x \in \Omega^c$ to
 a boundary point $q$, followed by a generalized billiard orbit
of length $t - s$  from $q$ back to $x$. By our assumption on the
length spectrum, the only possible orbit with length $t\in (r
L_{\gamma} - \epsilon, rL_{\gamma} + \epsilon)$ is $\gamma^r$ of
length $r$; but the  cutoff has removed this orbit. This completes
the proof.

\end{proof}

\subsection{Regularization}

In the preceding section we explained that the trace can be
microlocalized to $\gamma$. However, it is not true that the
operator $N$ can be microlcalized to $\gamma$, that is, we cannot
replace $N$ by $\chi_{\gamma} N \chi_{\gamma}$. This is due to the
singularities on the diagonal. We now explain how to regularize
them.

When we expand \begin{equation} \label{string} (N_{++ 0} + N_{++
1} + N_{+ -} + N_{-- 0} + N_{-- 1} + N_{-+})^M, \end{equation} we
get a sum of `strings' of factors. The  factors of  $N_{++0}$ or
$N_{-- 0}$ forbid the immediate application of the stationary
phase method due to the singularities of its phase and amplitude
on the diagonal. We now explain how to remove such factors and
regularize the term except for presence of $N_{++ 0}^M$ or $N_{--
0}^M$, which has to be handled separately. We will only explain
how to remove one factor of $N_{++ 0}$ next to a smooth factor
such as $N_{++ 1}$. In general, one just repeats the procedure one
pair at a time.

\subsection{The compositions $N_0 \circ N_1$}

We begin by characterizing the composition $N_0 \circ N_1$.
Essentially the same result holds for $N_1 \circ N_0$ by the same
argument.

The composed kernel equals
\begin{equation}\label{IKJ}  \begin{array}{l} N_0 \circ N_1 (k + i \tau, \phi_{1},  \phi_{2}) :
 =  (k + i \tau)^2 \int_{{\bf T}}
 \chi(k^{-1 + \delta} (\phi_{1} - \phi_{3} )) (1 - \chi(k^{-1 + \delta} (\phi_{2} - \phi_{3} ))) \\ \\
  H^{(1)}_1
 ((k \mu + i \tau) |q(\phi_{3 }) - q(\phi_{1})|)
\cos \angle (q(\phi_{3 }) - q(\phi_{1}), \nu_{q(\phi_{3 })} ) \\
\\
   H^{(1)}_1 ((k\mu  + i \tau) |q(\phi_{3} ) - q(\phi_{2} )|)
\cos \angle (q(\phi_{2}) - q(\phi_{3}), \nu_{q(\phi_{2})} ) d
\phi_{3} . \end{array}\end{equation}

The next proposition shows that these compositions are
semiclassical Fourier integral kernels. Thus, the singularity on
the diagonal is `regularized'. The crucial point is that the
composition with the singular factor $N_0$ decreased the order by
one unit.  We state in a somewhat simplified form and refer to
\cite{Z1} for a more detailed statement.

\begin{prop} \label{Onebadpairsym}  $N_0 \circ N_1 \circ \chi_0 (k + i \tau, \phi_{1},  \phi_{2}) )$
defines a semiclassical Fourier integral operator on $\partial
\Omega$ of order $-1$ associated to the billiard map.

%\end{itemize}

\end{prop}

\begin{proof}

We only sketch the proof and refer to \cite{Z1} for further
details.

The proof is an explicit calculation. Following \cite{AG}, we
change variables
$$\phi_3 \to \vartheta = \phi_1 - \phi_3,$$
 and then to change
variables
 $\vartheta \to u$, with:
\begin{equation}  u: = \left\{ \begin{array}{ll} |q(\phi_{3})
- q(\phi_{1})| , & \phi_{1} \geq \phi_{3} \\ & \\
- |q(\phi_{3}) - q(\phi_{1})|, & \phi_{1} \leq \phi_{3}
\end{array} \right.  = \left\{ \begin{array}{ll} |q(\phi_{1}
- \vartheta)
- q(\phi_{1})| , & \vartheta \geq 0 \\ & \\
- |q(\phi_{1} - \vartheta) - q(\phi_{1})|, & \vartheta \leq 0
\end{array} \right.\end{equation} In other words, we change from
the intrinsic distance along $\partial \Omega$ to chordal
distance.  The purpose of this change of variables is to simplify
the difficult factor in (\ref{IKJ}):

$$H^{(1)}_1 ((k\mu  + i \tau) |q(\phi_{3} ) - q(\phi_{1})|)
\to H^{(1)}_1 ((k\mu  + i \tau) |u|).$$

We then substitute the asymptotic WKB formula  for the `easy'
factor
 $H^{(1)}_1 ((k \mu + i \tau) |q(\phi_{1 }) - q(\phi_{3}|)$
 which is valid on the support of the cutoff.
The key points are then that

$$\left\{ \begin{array}{l}

(ii)  \cos (\angle q(\phi_{2}) - q( \phi_{2} + \vartheta),
\nu_{q(\phi_{2})})
 \to  |u| K(\phi_{2}, u), \;\; \mbox{with} \; K \; \mbox{smooth in } \; u;\\ \\
 (ii)
H^{(1)}_1 ((k \mu + i \tau)  |q(\phi_{1 }) - q(\phi_{2} -
\vartheta )|) \to      e^{i k |q(\phi_{1 }) - q(\phi_{2} )|}
  e^{i k  u a}
A(k + i \tau,  \phi_{ 1}, \phi_{2}, u), \\ \\
\mbox{where}\;\; A_k \; \mbox{is a symbol in } \;k\; \mbox{ of
order } \;\; -1/2 \;\;\mbox{ and smooth in }\; u. \end{array}
\right.
$$
Here,  \begin{equation} \label{a} a = \sin \vartheta_{1, 2},\;\;\;
\mbox{with} \;\; \vartheta_{1, 2} = \angle (q(\phi_{2}) -
q(\phi_{1}), \nu_{q(\phi_{2})}).
\end{equation}

 It follows that the composed kernel (\ref{IKJ})
 can be expressed in the form $A e^{i (k + i \tau) |q(\phi_1) -
 \phi_2)|}$ (further composed with $\chi_0$) with
\begin{equation} \label{AK} \begin{array}{l}
A(k + i \tau, \phi_{1}, \phi_{2})  = \int_{- \infty}^{\infty}
\tilde{\chi}(k^{1 - \delta} u) (1 - \chi(k^{1 - \delta}(\phi_2 -
\phi_1 - u)) \\ \\
\times  |u| e^{i k a u} H_1^{(1)}((k + i \tau) |u|) G((k + i
\tau), u, \phi_{1}, \phi_{2}) du,\end{array}
\end{equation} where $G((k + i \tau), u , \phi_{1}, \phi_{2})$ is a symbol
 in $k$ of order $-1/2$  and smooth in $u$,   and
where $\tilde{\chi}(k^{1 - \delta} u) = \chi(k^{1 - \delta}
(\phi_{1} - \phi_{3)}).$

We now  change variables again, $u' = k u$ (and then drop the
prime), to get
\begin{equation}\label{AKU} \begin{array}{lll} A(k + i \tau, \phi_{1},
\phi_{2})  & = & k^{-2}  \int_{- \infty}^{\infty}
\tilde{\chi}(k^{- \delta} u) (1 - \chi(k^{1 - \delta}(\phi_2 -
\phi_1 - k^{-1} u))) |u| e^{i a u} \\ & & \\
& \cdot&  H_1^{(1)}( b |u|) G(k + i \tau, \frac{u}{k}, \phi_{1},
\phi_{2}) du, \end{array} \end{equation} with $b = 1 + i(\tau/k)$.
   Since $|u| \leq k^{\delta}$ on the
support of the cutoff, we have $|\frac{u}{k} | \leq k^{-1 +
\delta}$. The Taylor expansion of $G(k + i \tau, u, \phi_{1},
\phi_{2})$ at $u = 0$ produces an asymptotic series  in $k$. The
main point is that the change of variables introduced a factor of
$k^{-2}$, which cancels the original factor of $k^2$ in
(\ref{IKJ}). We are left with an amplitude of order $-1/2$, which
defines a Fourier integral operator of order $-1$ with phase
$|q(\phi_1) -
 \phi_2)|$ as long as the integral (\ref{AKU}) has the same order
 as $G$.

After Taylor expanding, we end up with terms of the form
$$  \int_{- \infty}^{\infty} \tilde{\chi}(k^{- \delta} u)
|u| u^n  e^{i a u} H_1^{(1)}( b |u|)    du$$ times simple
functions of $(\phi_1, \phi_2)$.  The integral may be expressed in
the form
\begin{equation} \label{nthterm} \begin{array}{l}
\frac{\partial}{\partial b} \frac{\partial^n}{\partial a^n}
\int_{- \infty}^{\infty} \tilde{\chi}(k^{- \delta} u) e^{i  a u}
H_0^{(1)}(b |u|)    du |_{a = \sin \vartheta_{1, 2},  b = (1 +
i\tau/k)}. \end{array} \end{equation} Using the cosine transform
of the Hankel function, one can explicitly evaluate
(\ref{nthterm}) as $  \frac{\partial^n}{\partial a^n} (1 -
a^2)^{-3/2},$ modulo lower order terms. After composition with the
cutoff $\chi_0$ this factor is bounded above, so one gets that $A$
has order $-1/2$.

 \end{proof}

We iterate this proposition to deal with the full strings in
(\ref{string}). The only cases where this is not possible are the
ones which have no factors of $N_{+-}(k), N_{- +}(k)$. However,
they are still composed with the cutoff operator $\chi_{\gamma}$
and a somewhat similar calculation shows that this composition
produces a semiclassical pseudodifferential operator, whose order
is decreased from $0$  by one unit of $k$ for each   factor of
$N_0$,

\subsection{Final remarks and remainder estimate}

After removing the $N_0$ factors from the terms (\ref{string}) of
(\ref{GS}, we end up with a finite sum of ordinary oscillatory
integrals plus a complicated remainder.  In the next section, we
will describe the result of applying stationary phase to the
finite sum. We now make a few comments on the remainder, referring
to \cite{Z1} for further discussion.

The estimate of the remainder is complicated because the operator
norm $||N(k + i \tau)||$ does not decrease with increasing $\tau$.
To obtain a small remainder we redefine $\tau \to \tau \log k$.
This changes the wave trace expansions by  $k^{- C r L_{\gamma}}$
and hence for $R$ sufficiently large a remainder estimate of
$O(k^{-R})$ is sufficient to separate a finite part of the wave
trace from the remainder. We then estimate the remainder $$Tr
\int_{\R} \rho(k - \lambda)  N(\lambda)^{M_0 + 1}  (I \!+\!
N(\lambda \!+\! i \tau))^{-1} \; \chi_{\gamma} \; N'(k + i \tau) d
\lambda$$ by applying the Schwarz inequality for the
Hilbert-Schmidt inner product. We futher use standard  estimates
on the Poisson kernel to remove the factor of $(I \!+\! N(\lambda
\!+\! i \tau))^{-1}$. Unfortunately, this estimate will also
remove the cutoff $\rho$ and replaces $N^M$ by $N^M N^{*M}$. If
one regularizes these products as above, we find that the critical
points correspond to closed circuits of $M$-links which begin at
some point $x$, end at some point $x'$ and then return to $x$ by
traversing the links in reverse order. The cutoffs $\chi_0,
\chi_{\gamma}$  prohibit a proliferation of small links (grazing
rays) and  force the links in critical paths to point in the
direction of $\gamma$ and hence to be of length roughly $M
L_{\gamma}$. The imaginary part $i \tau \log k$ of the
semiclassical parameter then contributes a damping factor of $e^{-
\tau M L_{\gamma} \log k}$ for each link. The links correspond to
the $N_1$ factors. Thus, for each string,  we have one $k^{-1}$
for each $N_0$ factor and one  $e^{- \tau M L_{\gamma} \log k}$
for each $N_1$ factor.  For sufficiently large $\tau$ these
combine to give a factor of  $k^{-R}$ for any prescribed $R$ for
each term  of (\ref{string}).

\section{Calculating coefficients}
\label{COEFFS}

We now  apply the stationary phase method directly to each
regularized integrals. Additionally, we have to prove that the
remainder is as small as claimed.  The result is the following
explict formula for the wave invariants:

\begin{prop} We have:
$$\begin{array}{l}
\dis
B_{\gamma^r, j} + B_{\gamma^{-r}, j - 1}  = b_{r,j,2j} f^{2j}(0)
+ b_{r, j, 2j - 1} f^{2j - 1}(0) + \tilde{B}_{\gamma^r, j}, \;\;
\mbox{with}\\ \\
\qquad b_{r,j,2j} f^{2j}(0) + b_{r, j, 2j -1} f^{2j - 1}(0)  =
r \{ 2 (h^{11})^j f^{(2j)}(0)\\ [6pt]
\hspace{30pt}\dis +\ \{2  (h^{11})^j \frac{1}{2 \!-
\!2 \cos \alpha/2} +
 (h^{11})^{j \!-\! 2} \sum_{q = 1}^{2r} (h^{1 q})^3\} f^{(3)}(0) f^{(2j
\!-\! 1)}(0)\}\}\end{array}
$$ where $\tilde{B}_{\gamma^r, j}$ depends only on the $2j-2$
Taylor polynomial of $f$ at $x = 0.$ \end{prop}

We review  the calculation of the coefficients from  \cite{Z2} and
adapt it to the exterior case. We also take advantage of the
reduction to the boundary to simplify the arguments.

\subsection{Setting things up}

We recall that, in a small strip $T_{\epsilon}(\gamma)$ around
$\overline{a_+ a_-}$, the boundary $\partial \Omega$ consists of
two components which are symmetric graphs over the $x$-axis. We
write the graphs in the form $ y = \pm f(x)$ near $a_{\pm}$.

As in (\ref{MINT}), integrals over   $ (\Omega \cap T_{\epsilon}(\gamma))^M $consist of  $2^M$ terms, corresponding  to a choice of an element $\sigma$ of
$$\{\pm\}^M := \{\sigma : \Z_M \to \{\pm \} \}.$$
The length functional in Cartesian coordinates for a given assignment $\sigma$ of signs is given by
\setcounter{equation}{35}
\begin{equation} {\mathcal L}_{\sigma} (x_1, \dots, x_{2r}) = \sum_{j =
1}^{2r \!-\! 1} \sqrt{(x_{j \!+\! 1} \!-\! x_j)^2 \!+\!
(f_{\sigma(j \!+\! 1) }(x_{j \!+\!1}) \!-\!
f_{\sigma(j)}(x_j))^2}. \end{equation} If we write out the
integrals  $I_{M, \rho}^{\sigma}$ of (\ref{MINT2}) in Cartesian
coordinates, we obtain
\begin{equation} \label{GMG} \begin{array}{rl}
\dis\hspace{-15pt}
 I_{M, \rho}^{\sigma}(k  \!+\! i \tau)\!=\! \!\!\!\!\!&
\dis \int_{(-\epsilon, \epsilon)^M} \int_0^{\infty}
\int_{-\infty}^{\infty}  \hat{\rho}(t)
\bigg\{\Pi_{p = 1 }^{M - 1}  H^{(1)}_1 \big((k \lambda+ i \tau)\\[12pt]
&\times\ \dis |(x_{p \!+\! 1}
- x_{p}, f_{\sigma(p \!+\! 1)}(x_{p\! +\!1}) - f_{\sigma(p )}(x_{p}))|\big)
\\[6pt]
&\times\dis \frac{(x_{ p \!+\! 1} \!-\! x_{p}) \!-\! f_{\sigma(p \!+\! 1)}'(x_{p \!+
\!1})(
 f_{\sigma(p \!+\! 1)}(x_{p \!+\! 1}) \!-\! f_{\sigma(p )}(x_{p}))}{\sqrt{(x_{p \!+
\!1} \!-
\!x_{p})^2 \!+\! (f_{\sigma(p \!+\! 1)}(x_{p \!+\! 1}) \!+\! f_{\sigma(p )}(x_{p\!+\!1})(x_{p})
\!+\! L)^2}}    \bigg\}\\[12pt]
&\times\ \dis \chi(1 - \lambda)
e^{i k (1 - \lambda) t} d \lambda dt dx_1 \dots dx_{M}.
\end{array}
\hspace{-15pt}\end{equation}

 We now regularize the integrals
(\ref{GMG}) as above and obtain classical oscillatory integrals in
$x - y$ coordinates. We also
 eliminate the $(t, \lambda)$
variables by stationary phase. The Hessian in these variables is
easily seen to be non-degenerate, and the Hessian operator equals
$- \frac{\partial^2}{\partial t
\partial \lambda}.$ Since the amplitude depends on $t$ only in the
factor $\hat{\rho}(t),$ which is constant in a neighborhood of the
critical point,  only the zeroth order term in the Hessian
operator survives the method of stationary phase. We may therefore
eliminate the $ds d\mu$ integrals and replace the amplitude and
phase by their evaluations at $\mu = 1, t = r L_{\gamma}.$ This
simplifies the above to
\begin{equation}\label{EXPRESSION} \begin{array}{l}
k^{- \nu_{{\mathcal S}}}
 \!\!\!\int_{[-\epsilon, \epsilon]^{2r}}\!\!\!\!
e^{i k  {\mathcal L}_{\pm}(  x_1,  \dots, x_{2m})}
  \hat{\rho}({\mathcal L}_{\pm}(  x_1,  \dots, x_{2m})) a(k,   x_1,  x_2, \dots,  x_{2r})
d x_1 \cdots dx_{2r},\\ [6pt] . \end{array}
\hspace{-5pt}\end{equation}

An examination of the amplitudes and of the regularization
procedure leads to the following conclusions:
\begin{itemize}

\item (i)
The regularized integral $ I_{M, \rho}^{\sigma}((k + i \tau)) $ is
negligeable as $k \to \infty$  unless $M = 2 r$ (where $r
L_{\gamma}$ is the unique length in the support of $\hat{\rho}$),
and where only strings of $N_{- +}$ and $N_{+ -}$ are left upon
regularization. Otherwise there are no critical points.

\item (ii) The amplitudes of these oscillatory integrals have the
form
$$\sum_{n = 0}^{\infty} k^{-n} A_{ n}(x),$$
where $A_{ n}$ depends  only on the first $n + 2$ derivatives of
$f_{\sigma}.$

\item There is a {\it special   term}, namely the  initially regular  one in  which  $2r$
factors of $N_{-+}$ and $N_{+1}$ alternate, corresponding to $2r$
bounces of the bouncing ball orbit. The phase is simply ${\mathcal
L}_{\pm}$ and the amplitude is
\begin{equation}\label{AMPLSINGL}
\begin{array}{lll}
\dis a^0_{\pm}(k, x_1, \dots, x_{2r})
\ =\\[6pt]
\qquad \dis   \Pi_{p = 1}^{2 r }  a_1((k \!+\! i \tau) \mu
\sqrt{(x_{p } \!- \!x_{p\!+\!1})^2 \!+\! (f_{\sigma_{\pm}(p)}(x_{p
}) \!-\! f_{\sigma_{\pm}(p \!+ \!1)}(x_{p\!+\!1} )^2}) \\ [9pt]
\qquad \dis\times \ \frac{(x_{p} - x_{p+1}) f'_{\sigma_{\pm}(p)}
(x_{p}) -  ( f_{\sigma_{\pm}(p)} (x_{p}) -
f_{\sigma_{\pm}(p+1)}(x_{p+1}) )}{ \sqrt{(x_{p} - x_{p+1})^2 +
(f_{\sigma_{\pm}(p)} (x_{p}) - f_{\sigma_{\pm}(p+1)}(x_{p+1})
)^2}}\end{array} \end{equation} We only need the principal terms
in the symbols $a_1$, but to simplify notation we omit this step.

\end{itemize}

This amplitude and phase of the special term have certain  key
attributes which will be used extensively below

\eject\noindent
\begin{equation} \label{ATTRIBUTES}\left\{\! \begin{array}{l}
(i)~~~  \mbox{In  its dependence on}\; f, \; \mbox{ the amplitude
has the
form} \\[2pt]
\hspace{21pt}
{\mathcal A}(x, y, f, f'). \\  \\
(ii)~~  \mbox{In its dependence on the variables\ } x_p,
\mbox{the amplitude} \\[2pt]
\hspace{21pt} \mbox{has the form:}  \\ [4pt]
 \qquad a^0(k, x_1, \dots, x_{2r})   =   \Pi_{p = 1}^{2r} A_p (x_{p}, x_{p + 1})\;\;\;\; (2 p + 1 \equiv 1)\\  \\
(iii)~  D_{x_p}^{(2j - 1)} {\mathcal L}  \equiv  ((x_p - x_{p +
1})^2 +
 (f_{\epsilon_p}(x_p)\\[6pt]
\hspace{77pt} \dis -\  f_{\epsilon_{p+1}}(x_{p + 1}))^2)^{-1/2}
(f_{\epsilon_p}(x_p) \\[6pt]
\hspace{77pt} \dis -\ f_{\epsilon_{p+1}}(x_{p + 1}))
f_{\epsilon_p}^{(2j - 1)}(x_p)\\[6pt]
\hspace{74pt} \dis \implies   D_{x_p}^{(2j - 1)} {\mathcal L}
|_{x=0}
\equiv \epsilon_p f_{\epsilon_p}^{(2j - 1)}(0);\\ \\
(iv) \;\; \nabla a^0(k, x_1, \dots, x_{2r})_{x = 0} = 0.
\end{array} \right. \end{equation}
Here, the notation $\equiv$ means equivalence modulo terms with
fewer derivatives of $f_{\pm}.$ The last statement follows from
the values
\begin{equation} \label{VALUES} f(0) = L,  f'(0) = x|_{x = 0} = 0,  \angle(A - B, \nu_{B}) =
\pi/2, \end{equation} where $A, B$ denote the endpoints.

We now use this information to determine where the data
$f_{\pm}^{2j}(0),$\break $f^{(2j - 1)}_{\pm}$ first appears in the
stationary phase expansion for $I_k(a_{{\mathcal S}_0}, {\mathcal
L}_{\pm}) $.

\subsection{Stationary phase expansion}

 We now apply the method of stationary phase to these integrals
 and keep careful track of the dependence on the number of derivatives of
$f_{\pm}.$ To analyze the jungle of terms in the stationary phase
expansion, we use
 the diagrammatic approach (see e.g. \cite{AG}). We give a brief
 review.

Consider a general  oscillatory integral $Z_k = \int_{\Rr^n} a(x)
e^{ik S(x)} dx$ where $a \in C_0^{\infty}(\Rr^n)$ and where $S$
has a unique critical point in supp$a$  at $0$. Let us write $H$
for the Hessian of $S$ at $0$. The stationary phase expansion
takes the form:
$$\begin{array}{l}
\dis Z_k = \bigg(\frac{2\pi}{k}\bigg)^{n/2} \frac{e^{i \pi sgn
(H)/4}}{\sqrt{|det
H|}} e^{i k S(0)} Z_k^{h \ell}, \\ [6pt]
\hspace{70pt}\mbox{where}\;\;
\dis Z_k^{h \ell} =  \sum_{j = 0}^{\infty} k^{-j}
\bigg\{\sum_{(\Gamma, \ell):
\chi_{\Gamma'} = j} \frac{I_{\ell} (\Gamma)}{S(\Gamma)}\bigg\}.
\end{array} $$

Here, the sum runs over the  set
${\mathcal G}_{V, I}$  of  labelled graphs $(\Gamma, \ell)$ with $V$ closed  vertices of valency $\geq 3$ (each corresponding
to the phase), with one open vertex (corresponding to the amplitude), and with $I$ edges.
Further, the  graph $\Gamma'$ is defined to be $\Gamma$ minus the open vertex,
and $\chi_{\Gamma'} = V - I$ equals its Euler characteristic.
We note that there are only finitely many graphs for each $\chi$ because the valency condition forces
$I \geq 3/2 V.$ Thus, $V \leq 2 j, I \leq 3 j.$

The function $\ell$ `labels' each end
of each edge of $\Gamma$ with an index $j \in \{1, \dots, n\}.$
Also, $S(\Gamma)$ denotes  the order of the automorphism group of $\Gamma$, and  $I_{\ell} (\Gamma)$ denotes the `Feynman
amplitude' associated to $(\Gamma, \ell)$. By definition, $I_{\ell}(\Gamma)$  is obtained by the following rule: To each edge
with end labels $j,k$ one assigns a factor of $\frac{-1}{ik} h^{jk}$ where $H^{-1} = (h^{jk}).$ To each closed
vertex one assigns a factor of $i k \frac{\partial^{\nu} S (0)}{\partial x^{i_1} \cdots \partial x^{i_{\nu}}}$ where
$\nu$ is the valency of the vertex and $i_1 \dots, i_{\nu}$ at the index lables of the edge ends incident on
the vertex. To the open vertex, one assigns  the factor  $\frac{\partial^{\nu} a(0)}{\partial x^{i_1} \dots \partial x^{i_{\nu}}}$, where $\nu$ is its valence.   Then $I_{\ell}(\Gamma)$
is the product of all these factors.  To the empty
graph one assigns the amplitude $1$.  In summing over $(\Gamma, \ell)$ with a fixed
graph $\Gamma$, one sums the product of all the factors as the indices run over
$\{1, \dots, n\}$.

\subsubsection{The special term: The data
\boldmath{$f_{\pm}^{2j}(0)$}}

We first claim that  $f_{\pm}^{(2j)}(0)$ appears first in the
$k^{-j + 1}$ term. This is because any labelled  graph $(\Gamma,
\ell) $ for which $I_{\ell}(\Gamma)$  contains the factor
$f_{\pm}^{(2j)}(0)$ must have a closed vertex of valency $\geq
2j$,  or the open vertex must have  valency $\geq 2j - 1.$ The
minimal absolute Euler characteristic $|\chi(\Gamma')|$ in the
first case is $1 - j$. Since the Euler characteristic is
calculated after the open vertex is removed, the minimal absolute
Euler characteristic in the second case is $-j$ (there must be at
least $j$ edges.) Hence such graphs do not have minimal absolute
Euler characteristic. It follows that the  only  labelled graph
$(\Gamma, \ell)$ with $-\chi(\Gamma') = j - 1$ with $I_{\ell}
(\Gamma)$ containing  $f_{\pm}^{(2j)}(0)$ is  given by:
\begin{itemize}

\item ${\mathcal G}_{1, j}^{2j, 0} \subset {\mathcal G}_{1, j} $with  $V = 1, I = j;$
the unique graph has no open vertex, one closed vertex and    $j$ loops at the  closed vertex.
 The only labels producing the desired data are those $\ell_p$ which assign
all  endpoints
of all edges labelled  the same index $p$.

\item Assuming the up/down symmetry, the corresponding Feynman amplitude $I_{\ell_p}(\Gamma)$  has the form
$ 2 r L   (h^{11})^j  f^{(2j)}(0) + \cdots$,
where again $\cdots$ refers to terms with $\leq 2j - 1$ derivatives.

\end{itemize}

To prove the last statement about the Feynman amplitude,
we  first observe that
$f^{(2j)}(0)_{\pm}$ appears
linearly in $I_{\ell_p}(\Gamma)$.
Indeed,  the only labelled graphs which produce this datum give the same   label to all endpoints of
all edges,   corresponding to applying only
 derivatives in a single variable $\frac{\partial}{\partial x_k}$.
 Using (\ref{VALUES}),  an examination of (\ref{AMPLSINGL}) shows
that
 the coefficient of $f^{(2j)}(0)$ of $I_k(a_{{\mathcal S}_0}, {\mathcal L}_+)$ equals
$$\begin{array}{l}
\dis\sum_{p = 1, p \equiv 1}^{2r }  (h^{pp}_+)^j
\bigg(\frac{\partial}{\partial x_p}
\bigg)^{2j}  {\mathcal L}_+ ( \dots, x_k, \dots x_{2r})
= L
\bigg[\sum_{p = 1, p \equiv 1}^{2r }  (h^{pp}_+)^j
\bigg]  f^{(2j)}_{+}(0). \end{array} $$
To complete the proof it suffices to  observe that   the diagonal matrix elements $h^{pp}_{\pm}$ are constant
in $p$ and in the  sign  $\pm.$

Had we assumed that the obstacles were both up-down and right-left
symmetric as in \cite{Z3}, then we would already  have solved the
inverse spectral problem in this step,  since all odd Taylor
coefficients of $f$ vanish at $x = 0$ and we will be able to
recover all even ones.

\subsubsection{The special term: The data
\boldmath{$f_{\pm}^{(2j - 1)}(0)$}}

We now consider the more difficult odd coefficients  $f_{\pm}^{(2j
- 1)}(0)$, which will require the attributes of
 the amplitude (\ref{AMPLSINGL}) detailed in (\ref{ATTRIBUTES}).

We again claim that the Taylor coefficients $f_{\pm}^{(2j -
1)}(0)$ appear first in the term of order $k^{-j + 1},$ and  we
enumerate the labelled graphs which have a combinatorial structure
capable of producing $f_{\pm}^{(2j - 1)}(0)$ as a factor in
$I_{\ell}(\Gamma)$ in the  $k^{-j + 1}$  term. In fact, only the
main term  $I_{M, \rho, C}^{\sigma_{\pm}}$ will produce such a
term. The following is proved in \cite{Z2}:

\begin{lem} We have:

\noindent{\bf (i)} There are no labelled graphs with $-\chi':= - \chi(\Gamma') \leq j - 1$ for which $I_{\ell}(\Gamma)$ contains
the factor $f_{\pm}^{(2j - 1)}(0)$.  \\
\noindent{\bf (ii)} There are exactly two types of labelled diagrams $(\Gamma, \ell)$ with
 $\chi(\Gamma') = -j + 1$ such that  $I_{\ell}(\Gamma)$ contains the factor.
  $f_{\pm}^{(2j - 1)}(0)$. They are given by:

\begin{itemize}

\item ${\mathcal G}_{2, j + 1 }^{2j - 1, 3, 0} \subset {\mathcal G}_{2, j + 1 }$ with $V = 2, I = j + 1$):
Two closed vertices, $j - 1$ loops at one closed  vertex, $1$ loop at the second closed vertex, one edge between the closed vertices; no open vertex.   Labels $\ell_{p,q}$: All labels at the closed vertex with valency $2j - 1$ must be the same
index $p$ and all at the closed vertex must the be same index $q$. Feynman amplitude:
$(h^{pp}_{\pm})^{j - 1} h^{qq}_{\pm} h^{pq}_{\pm}  D_{x_p}^{2j - 1} {\mathcal L}_{\pm} D_{x_q}^3 {\mathcal L}_{\pm}
\sim (h^{pp}_{\pm})^{j - 1} h^{qq}_{\pm} h^{pq}_{\pm} f^{(2j - 1)}(0) f^{(3)}(0).  $

\item ${\mathcal G}_{2, j + 1 }^{2j - 1, 3, 0} \subset {\mathcal G}_{2, j + 1 } $ with $ V = 2, I = j + 1$:
Two closed vertices, with $j - 2$ loops at one closed vertex, and with three edges between the two closed vertices; no open vertex.   Labels $\ell_{p,q}$: All labels at the closed vertex with valency $2j - 1$ must be the same
index $p$ and all at the closed vertex must the be same index $q$; $(h^{pp}_{\pm})^{j - 2} (h^{pq}_{\pm})^3 D_{x_p}^{2j - 1} {\mathcal L}_{\pm} D_{x_q}^3 {\mathcal L}_{\pm} \sim (h^{pp}_{\pm})^{j - 2} (h^{pq}_{\pm})^3 f^{(2j - 1)}(0) f^{(3)}(0). $

\end{itemize}

\end{lem}

This Lemma requires a careful  search through the diagrams  and an
explicit calculation of the amplitudes. It is the most `unstable'
step of the proof, in which any calculational error could ruin the
result.  For instance, it turns out that in addition to the two
diagrams which contribute, there are three more which have a
combinatorial structure capable of contributing the datum $f^{2j -
1}(0)$ ; but the corresponding amplitudes turn out to vanish due
to special attributes (\ref{ATTRIBUTES}) of the Taylor expansion
of the amplitude at the critical point. We should note other
attributes in  addition to (\ref{ATTRIBUTES}) are needed to show
that the Feynman amplitude for these diagrams vanishes.

\subsubsection{Other terms}
We also have:
\begin{lem} Other terms do not contribute the data
 $f^{(2j)} (0), f^{(2j -
1)}(0)$  in the term of order $k^{-j + 1}.$ \end{lem}

The proof of this is easy. It merely requires checking how the
regularization compares to the stationary phase expansion. As
noted above, the regularized integrals are lower order yet do not
involve higher derivatives of the amplitude or phase.

\section{Recovering the domain}
\label{DOMAIN}

In this section, we sketch the proof of the main result, Theorem
\ref{ONESYM}.  The proof is virtually identical to the proof that
one can determine a $\Z_2$-symmetric analytic domain from the wave
invariants of a $\Z_2$- symmetric hyperbolic bouncing ball orbit.
In fact, the proof requires fewer assumptions, since
 $L_{\gamma}$ is a resonance invariant (being  the least non-zero singular time in the wave trace.)
We may then inspect the wave invariants (i.e. the Balian-Bloch trace invariants) at the iterates $\gamma^r$ of $\gamma.$  Give the expressions above for the wave invariants,
it suffices to analyse the Hessian coefficients.

\subsection{Poincare map and Hessian of the length functional at
\boldmath{$\gamma$}}

We will need to recollect certain facts about the Hessian of the
length function at critical points corresponding to the iterates
$\gamma^r$. In most respects,the discussion is   identical to the
interior case, so we refer the reader  to \cite{Z1, Z2} for
background on the linear Poincare map $P_{\gamma}$ and its
relation to the Hessian of the length function. We will only
consider the $\Z_2$-symmetric case since it is the one relevant to
our result.

Henceforth we assume ${\mathcal O}$ is convex,  so that the
bouncing $\gamma$ is necessarily hyperbolic and  the eigenvalues
of its Poincare map $P_{\gamma}$ are of the form $\{e^{\pm
\lambda}\}.$ They are related to the geometry by
\setcounter{equation}{42}
\begin{equation} \label{COSHALPHA}
\cosh \alpha/2 = \bigg(1 - \frac{L}{R}\bigg), \end{equation}
where $R$ is the common radius  of curvature at the endpoints
$a_j$. For details, see  \cite{PS}\cite{KT}. When ${\mathcal O}$
is not convex, $\gamma$ could be elliptic. We refer to \cite{Z2}
for the calculation of the wave invariants in that case.

  We will need formulae involving matrix  elements of the inverse Hessian matrix $H_{2r}^{-1} = (h^{pq})$
at the $\Z_2$-symmetric bouncing ball orbit.
As in \cite{Z2}, the Hessian in $x-y$ coordinates is the {\it circulant } matrix
\begin{equation} H_{2r} = C(2 \cosh \alpha/2, 1, 0, \dots, 0, 1), \end{equation}
that is, the matrix whose  rows (or columns) are obtained by
cycling the first row (or column) (see \cite{Z1, Z2} for
background). Therefore $H_{2r}$ is diagonalized by the finite
Fourier matrix $F$ of rank $2r$ (see \cite{Z1, Z2} for its
definition and background). Exactly as in \cite{Z1, Z2}, we have
\begin{equation} \label{CIRCHESS} H_{2r} = F^* \mbox{diag}\;
\bigg(2 \cos \alpha/2 + 2 ,
\dots, 2 \cos \alpha/2 + 2 \cos \frac{(2r - 1) \pi}{r}\bigg) F.
\end{equation}

An important role in the inverse problem is played by sums of powers of elements of columns of $H_{2r}^{-1}.$ From
(\ref{CIRCHESS}), it follows that
\begin{equation} \label{INVHESS} \begin{array}{l}
\dis H_{2r}^{-1}
= F^* \bigg(\mbox{diag}\;
\bigg(\frac{1}{2 \cos \alpha/2 + 2} , \dots,\frac{1}{
2 \cos \alpha/2 + 2 \cos \frac{(2r - 1) \pi}{r}} \bigg)\bigg) F.
\end{array}
\hspace{-10pt}\end{equation}
Hence,  the first row $[H_{2r}^{-1}]_1 = (h^{11}, \dots, h^{1
2r})$ (or column) of the inverse is given by:
\begin{equation} \label{ROW} h^{1 q} = \sum_{k = 0}^{2r - 1} \frac{w^{ (q - 1) k } }{p_{\alpha, r}(w^k)}, \end{equation}
where $w = e^{\frac{2 \pi i }{2r}}$ and where $p_{\alpha, r} (z) =
2 \cosh \alpha/2 +  z  +  z^{n-1}$, in \cite{Z1, Z2} we will need
that, for any $p$,
\begin{equation}\label{FIRSTPOWER} \sum_{q = 1}^{2 r }  h^{pq} = \frac{1}{ 2 +  \cosh \alpha/2}.   \end{equation}

\subsection{Domain recovery}
We now prove by induction on $j$ that the Taylor coefficients
$f^{2j-1}(0), f^{2j}(0)$ can be determined from $B_{\gamma^r,j}$
as $r$ varies over $r = 1, 2, 3, \dots.$ At first, we assume
$f^{3}(0) \not= 0.$

It suffices to  separately determine the two terms
\begin{equation} \label{TWOTERMS}
\ba{lll}\dis
2 (h^{11}_{2r})^2
\bigg\{f^{(2j)}(0) + \frac{1}{2 - 2 \cos \alpha/2}\ f^{(3)}(0)
f^{(2j - 1)}(0) \bigg\}, \\[9pt]
\dis
\mbox{and} \quad
\bigg\{  \sum_{q = 1}^{2r} (h^{1 q}_{2r})^3
\bigg\} f^{(3)}(0) f^{(2j -
1)}(0).
\ea
\end{equation}
As discussed in
\cite{Z2}, the terms decouple   as $r$ varies if and only if\break
 $F_r(\cos \alpha/2) : =  \sum_{q = 1}^{2r} (h^{1 q}_{2r})^3$ is non-constant in $r = 1, 2, 3, \dots$ .

By the explicit calculation in \cite{Z2}, we  have:
$$
\ba{clll}\dis
\sum_{q =
1}^{2 r } ( h^{pq})^3 =\\[9pt]
\dis 2r  \!\!\!\sum_{k_1, k_2 = 0}^{2r} \frac{1
}{(\cosh \alpha/2 \!+\! \cos \frac{k_1 \pi}{r}) (\cosh \alpha/2 \!+ \!\cos
\frac{k_2 \pi}{r})(\cos \alpha/2 \!+\! \cosh \frac{(k_1 \!+\!k_2)
\pi}{r})} .
\ea$$
It is obvious that the sum is strictly increasing as
$r$ varies over even integers.

      We now begin the inductive argument. From the  $j = 0$ term we determine $f''(0).$ Indeed,
 $(1 - L f^{(2)}(0) = \cos(h) \alpha/2$ and $\alpha$ is a wave trace invariant.
  From the $j = 2$ term we recover $f^{3} (0), f^{4}(0).$  The induction hypothesis is then that the Taylor polynomial
  of $f$ of degree $2j - 2$ has been recovered by the $j - 1$st stage. By the decoupling argument we can determine
  $f^{2j}(0), f^{2j - 1}(0)$ as long as $f^{3}(0) \not= 0.$

   As discussed in \cite{Z2}, the domain can still be recovered
  if $f^{3}(0) = 0,$ but $f^5(0) \not= 0$. The argument is
  essentially the same. Similarly if $f^{3}(0) = f^5(0) = 0,$ but $f^7(0) \not=
  0$, and so on. If all odd derivatives vanish, then one still
  recovers the domain as noted above.

This completes the sketch of the  proof of Theorem (\ref{ONESYM}).

\end{document}